\documentclass[a4paper, 11pt]{article}

\usepackage{amssymb,latexsym,amsmath}     
\usepackage{graphicx,color,afterpage,authblk}

\usepackage{algorithm,algpseudocode} 
\usepackage{mathtools}
\DeclarePairedDelimiter\ceil{\lceil}{\rceil}

\newtheorem{remark}{Remark}

\newcommand{\R}{\mathbb{R}}

\graphicspath{{./}{figures/}}

\begin{document}

\title{Efficient implementation of characteristic-based schemes \\ on unstructured triangular grids}

\author[1]{S. Cacace \thanks{\tt cacace@mat.uniroma3.it}}
\author[1]{R. Ferretti \thanks{\tt ferretti@mat.uniroma3.it}}
\affil[1]{Dipartimento di Matematica e Fisica, Universit\`a Roma Tre, Largo S. Leonardo Murialdo, 1, Roma, Italy}

\maketitle

\begin{abstract}
Using characteristics to treat advection terms in time-dependent PDEs leads to a class of schemes, e.g., semi-Lagrangian and Lagrange--Galerkin schemes, which preserve stability under large Courant numbers, and may therefore be appealing in many practical situations. Unfortunately, the need of locating the feet of characteristics may cause a serious drop of efficiency in the case of unstructured space grids, and thus prevent the use of large time-step schemes on complex geometries.

In this paper, we perform an in-depth analysis of the main recipes available for characteristic location, and propose a technique to improve the efficiency of this phase, using additional information related to the advecting vector field. This results in a clear improvement of execution times in the unstructured case, thus extending the range of applicability of large time-step schemes.

\end{abstract}

\noindent{\bf Keywords:}
Large time-step schemes, unstructured grids, point location, computational complexity.\\

\noindent{\bf AMS subject classification:} 65-04, 65D18, 65M06, 65M25.

\section{Introduction}

Born in the 50s in the framework of environmental fluid dynamics and Numerical Weather Prediction, large time-step, characteristic-based schemes have become in recent years a useful tool for various PDE models, mainly of hyperbolic type. While this class of schemes collects various techniques (for example semi-Lagrangian \cite{FF13}, Lagrange--Galerkin \cite{DR82,P82}, ELLAM \cite{RC02}) having in common the use of the method of characteristics to treat advection terms, to fix ideas we will refer in what follows to the case of semi-Lagrangian (SL) schemes, which probably employ this strategy in its simplest form.
We consider, as a model problem, the simple variable-coefficient advection equation,
\begin{equation}\label{eq:trasp}
\begin{cases}
u_t + f(x,t)\cdot \nabla u = 0 & (x,t)\in\R^d\times\R^+, \\
u(x,0) = u_0(x),
\end{cases}
\end{equation}
in which we choose $d=2$ (the general case is conceptually similar). The solution of \eqref{eq:trasp} may be represented via the well-known formula of characteristics
\begin{equation}\label{eq:caratt}
u(x,t) = u_0(X(x,t;0)),
\end{equation}
where $X(x,t;s)$ is the solution at time $s$ of the ordinary differential equation
\[
\begin{cases}
\frac{d}{ds} X(x,t;s) = f(X(x,t;s),s) & s\in\R, \\
X(x,t;t) = x,
\end{cases}
\]
that is, the trajectory moving with velocity $f(X,s)$ and passing through the point $x$ at time $t$. We assume that $f$ is $C^1$ with bounded derivatives on the whole of $\R^2$, so that:
$$
\|f(x_1,t_1)-f(x_2,t_2)\| \le L_x \|x_1-x_2\|+L_t \|t_1-t_2\|,
$$
with $L_x$, $L_t$ denoting the two Lipschitz constants associated to respectively space and time increments, and $\|\cdot\|$ denoting the Euclidean norm. Clearly, such a framework is ultimately directed towards nonlinear equations in which the advection term has a smooth space and time dependence, at least in a large majority of the computational domain.

Once a time grid $t_n=n\Delta t$ has been set, a SL discretization of \eqref{eq:trasp} uses the representation formula \eqref{eq:caratt} written on a single time step, i.e.,
\begin{equation*}
u(x,t_{n+1}) = u(X(x,t_{n+1};t_n),t_n).
\end{equation*}
To turn this relationship into a computable scheme, we build a space grid with space scale $\Delta x$ and with nodes in the set $\mathcal V=\{x_i\}_{i=1,\dots,N}$, the foot of the characteristic $X(x_i,t_{n+1};t_n)$ is replaced by a numerical (e.g., one-step) approximation $X^\Delta(x_i,t_{n+1};t_n)$, and the value $u(\cdot,t_n)$ by an interpolation $I[V^n](\cdot)$, constructed using the vector $V^n=(v_1^n \cdots v_N^n)$ of the node values at time $t_n$, with $v_i^n$ corresponding to the $i$-th node $x_i$ and the $n$-th time step $t_n$, where $N=|\mathcal V|$ denotes the total number of nodes.

The scheme is therefore in the form
\begin{equation}\label{eq:schema_trasp}
v_i^{n+1} = I[V^n]\left(X^\Delta(x_i,t_{n+1};t_n)\right).
\end{equation}
In \eqref{eq:schema_trasp}, the discrete approximation $X^\Delta(x_i,t_{n+1};t_n)$ of $X(x_i,t_{n+1};t_n)$ might be obtained in the simplest case by applying, backward in time, the explicit Euler scheme:
\begin{equation}\label{eq:eulero}
X^\Delta(x_i,t_{n+1};t_n) = x_i-\Delta t \> f(x_i,t_{n+1}).
\end{equation}
When working on a bounded domain  $\Omega\subset\R^2$, boundary conditions (e.g. of  Dirichlet type) can be treated by a suitable variable step technique, which basically stops characteristics on the boundary, as discussed in \cite{FF13}. This only modifies the definition of the points $X^\Delta(x_i,t_{n+1};t_n)$, and therefore we will not give further details here. 

Concerning the interpolation, this step is typically accomplished in local form, using the values of the numerical solution at nodes close to $X^\Delta(x_i,t_{n+1};t_n)$. Selecting the relevant values requires an $\mathcal{O}(1)$ cost on a structured array of nodes, and therefore is not a critical issue from the viewpoint of complexity. On the other hand, when working on unstructured (typically, but not necessarily, triangular) grids, the interpolation is usually computed via Lagrange finite elements: interpolating at a given point requires to first select the element containing the point, and then use the Lagrange basis associated to this specific triangle. In comparison with the structured case, in the unstructured case the former phase (point location) represents a clear bottleneck, which either prevents the use of large time-step schemes, or causes a substantial drop in their efficiency. In fact, as we will show in the last section, the point location phase covers a significant part of the total CPU time.

Despite this difficulty, a certain amount of literature has been devoted to unstructured implementations of characteristic-based schemes; in most cases, however, we found that an in-depth discussion of the efficiency issues is eluded. In other cases, practical recipes are provided: the two typical techniques used are on one hand the quadtree search (see \cite{G98,G00}), on the other the tracking of characteristics via substepping, which requires in general to move from one element to its neighbour, thus making the search easier (see \cite{RBS06} for the case of a triangular mesh, \cite{BDR13,B20} for a Voronoi mesh). We will briefly review the ideas behind these techniques in the next section.

To the authors' knowledge, the optimal complexity of known general-purpose point location algorithms is $\mathcal{O}(\log N)$, where $N$ is the number of grid nodes. In this paper, we will show that this complexity may be brought to $\mathcal{O}(1)$, by using the information related to the specific problem under consideration, that is, moving from a general-purpose algorithm to an algorithm tuned on the case of characteristics, at the price of introducing some additional data structures related to the mesh. 
A first motivation for this study is to apply efficient semi-Lagrangian techniques to Navier--Stokes equations on non-orthogonal geometries \cite{BCCF18}.

The paper is structured as follows. In Section \ref{sec:location}, we review the two main techniques to locate a point in an unstructured triangulation, and study their computational complexity. In Section \ref{sec:new}, we study in detail some possibilities to improve the point location algorithms. Last, in Sections \ref{sec:numeric} and \ref{sec:conclu} we present a numerical validation for the algorithm and draw some conclusions.

\section{Locating a point on a triangular grid: some basic facts}\label{sec:location}

In this section, we briefly review two major approaches to point location on triangulations, namely the quadtree search and the barycentric walk search, including an experimental analysis of their computational complexities.

\subsection{Quadtree algorithm}\label{QT-rules} As far as the authors know, the first appearance of this algorithm dates back to the 70s \cite{FB74}. The algorithm is based on an auxiliary data structure of quadtree type, i.e., a tree where all nodes but the leaves have precisely four children. Each node (also termed as a {\em quad}) corresponds to a rectangle, starting with the quad associated to the root and containing the whole triangulation, and each successive level divides the quad into four. Once fixed an integer $q\ge 2$, the subdivision is stopped as soon as one of the following conditions is satisfied (see Fig. \ref{fig:quadtree}):

	\begin{figure}
	\centering
	\includegraphics[height=6cm]{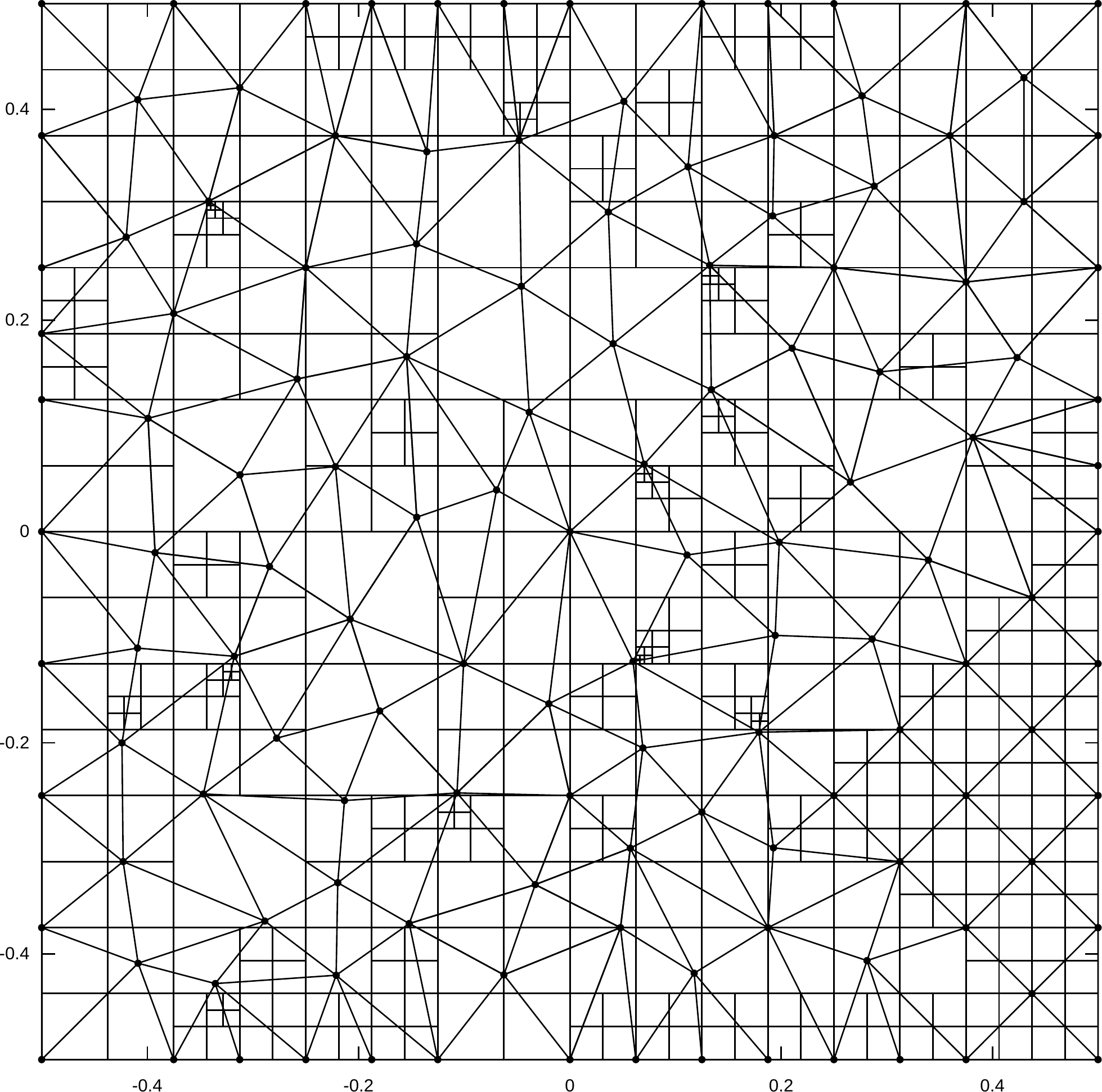}
	\caption{Quadtree partition associated with an unstructured triangulation of $\Omega=[-\frac12,\frac12]^2$, with $q=3$.}\label{fig:quadtree}
	\end{figure}

\begin{enumerate}

\item[$a)$]
The quad intersects a number $n_t$ of triangles such that $1\le n_t\le q$, and contains no vertex;

\item[$b)$]
The quad contains exactly one vertex, regardless of the number of triangles $n_v$ it joins;

\item[$c)$]
The quad does not intersect the triangulation.

\end{enumerate}

A leaf of the tree is generated at the final level of the subdivision, and the list of triangles intersecting the final quad is associated to the leaf.
A point location requires to visit the tree: once found the leaf containing the point, the location of the point in the triangulation is completed with a number $\mathcal{O}(\max(q,n_v))$ of operations. The tree is unbalanced in general; however, for a regular Delaunay triangulation we can reasonably assume that the average complexity for the visit of the tree (and therefore, for one point location) is $\mathcal{O}(\log N)$, where $N$ is the number of grid nodes, while the complexity of the checks to be done at the leaves is constant. The complexity of a single point location takes then the form
$$
\mathcal{O}\left(C_1^Q + C_2^Q \log N\right).
$$
In principle, complexity of the visit should depend on $q$; however, a decrease of $q$ causes at the same time a higher depth of the quadtree and a shorter list of elements to be checked at a leaf, and vice versa for an increase of $q$ (for example, in the grid of Fig. \ref{fig:quadtree}, the relatively low value of $q=3$ causes a tree depth of ten levels with only 218 elements). Except for the lowest values of $q$, which may lead to an extremely deep tree, the two effects tend to compensate, as shown by the following numerical test.

\begin{figure}[!h]
	\centering
	\includegraphics[width=0.47\textwidth]{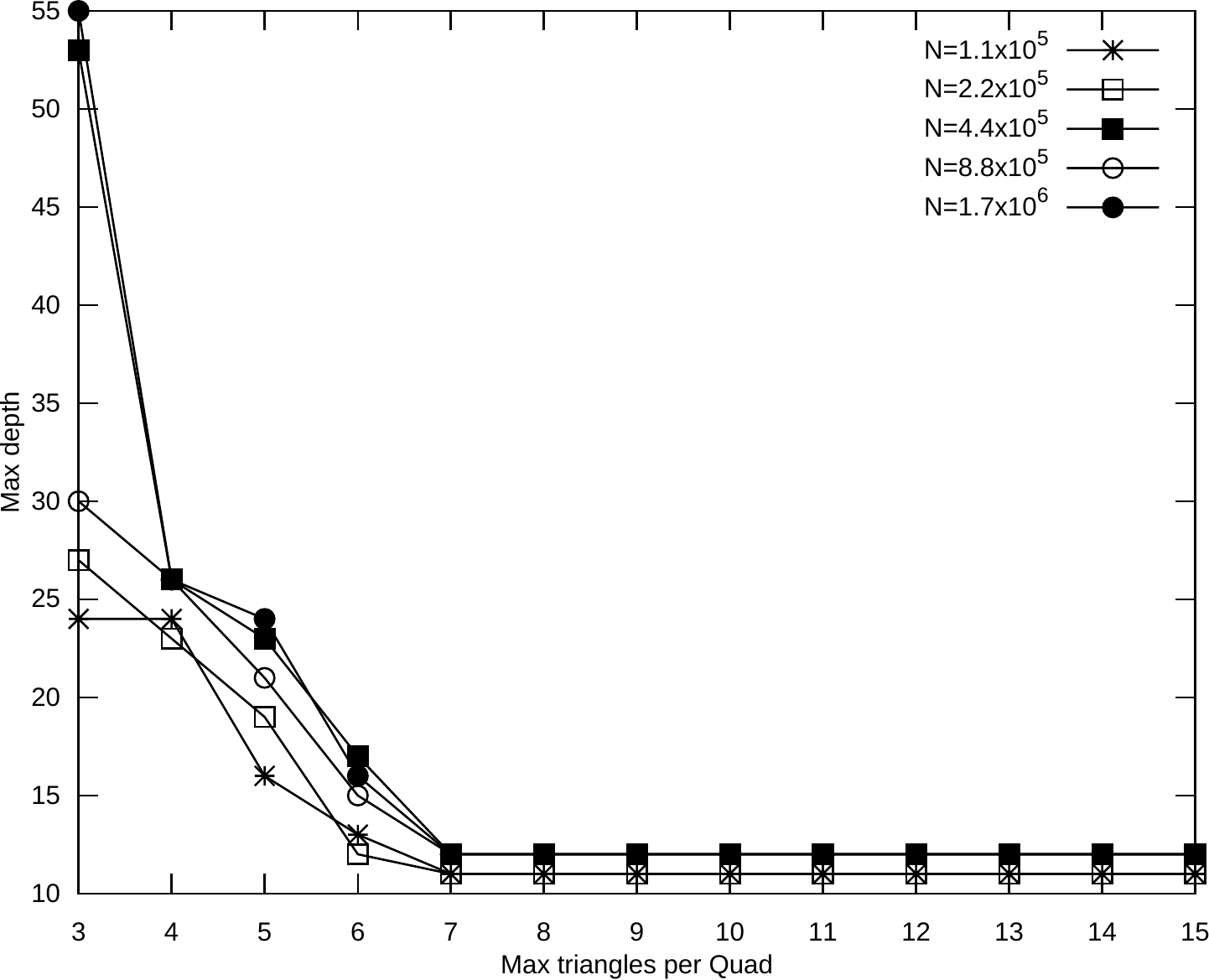} \quad
	\includegraphics[width=0.47\textwidth]{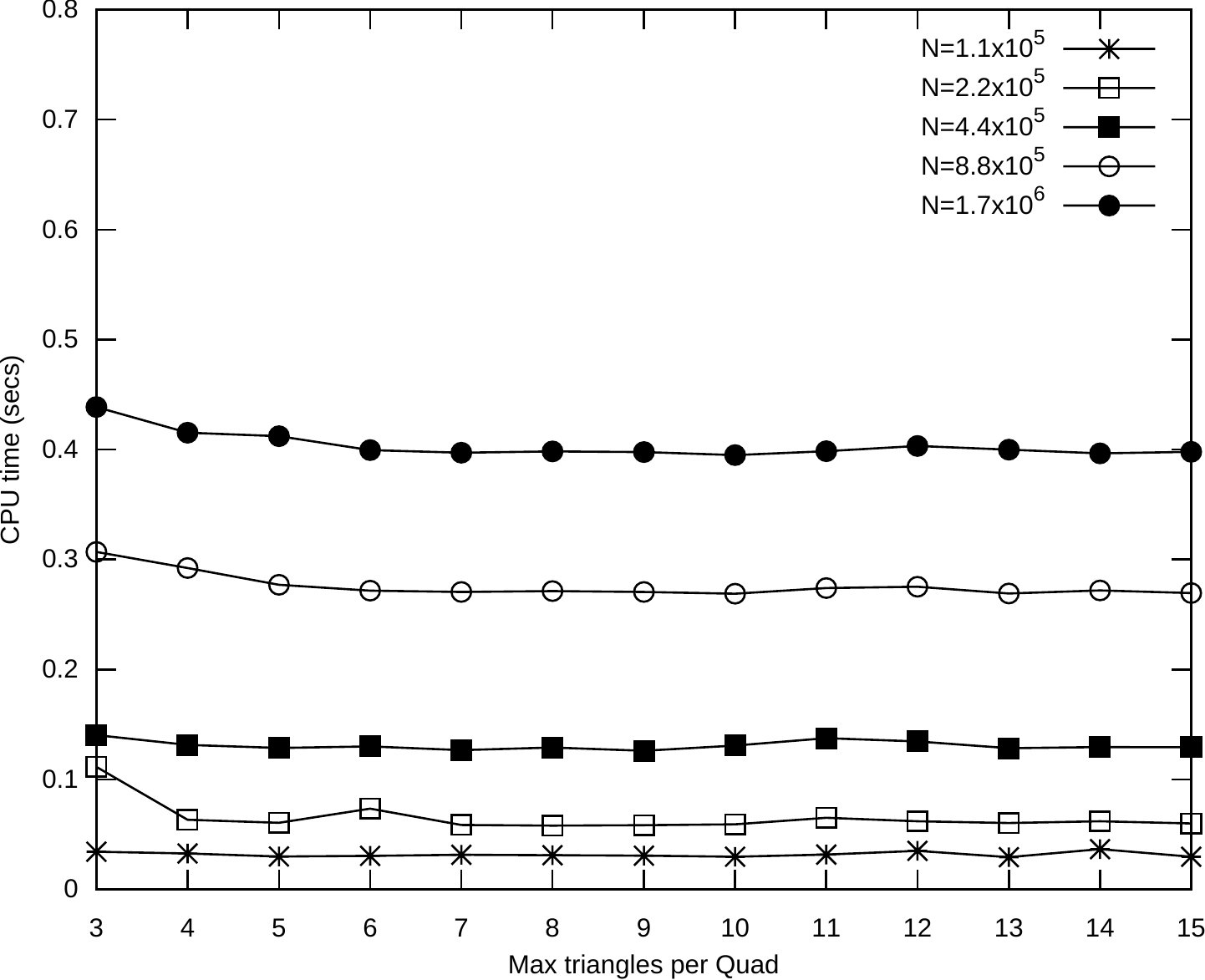}
	\caption{Depth of the quadtree versus $q$ (left) and search times versus $q$ (right), for $10^5\lesssim N\lesssim 1.7\cdot 10^6$.}\label{fig:qt_q}
	\end{figure}

\paragraph{Quadtree: numerical example.} We show here an experimental assessment of the performance of quadtree search. In the first plots (Fig. \ref{fig:qt_q}), we consider meshes of size $N$ ranging from about $10^5$ to about $1.7\cdot 10^6$. In the left plot, we compare the depth of the tree obtained for different values of $q$: the plot shows a clear saturation effect, and for $q\ge 7$ the depth becomes constant for all meshes. This effect might be explained with the fact that, in a regular Delaunay mesh, this is the typical maximum number of triangles joined at a node, this meaning that leaves of both types $(a)$ and $(b)$ contain typically a similar number of triangles (no more than $q$). Then, the tree is likely to be more balanced, and at the increase of $q$ we don't expect any improvement in the depth of the tree. 

In the right plot, we compare the search time for a set of $N$ random query points (the same size of the mesh), obtained with different values of $q$; here, it is clear that this parameter has a small effect, if any at all, on the execution times. From now on, we will choose $q=7$ in all the tests. 

\begin{figure}
	\centering
	\includegraphics[height=6cm]{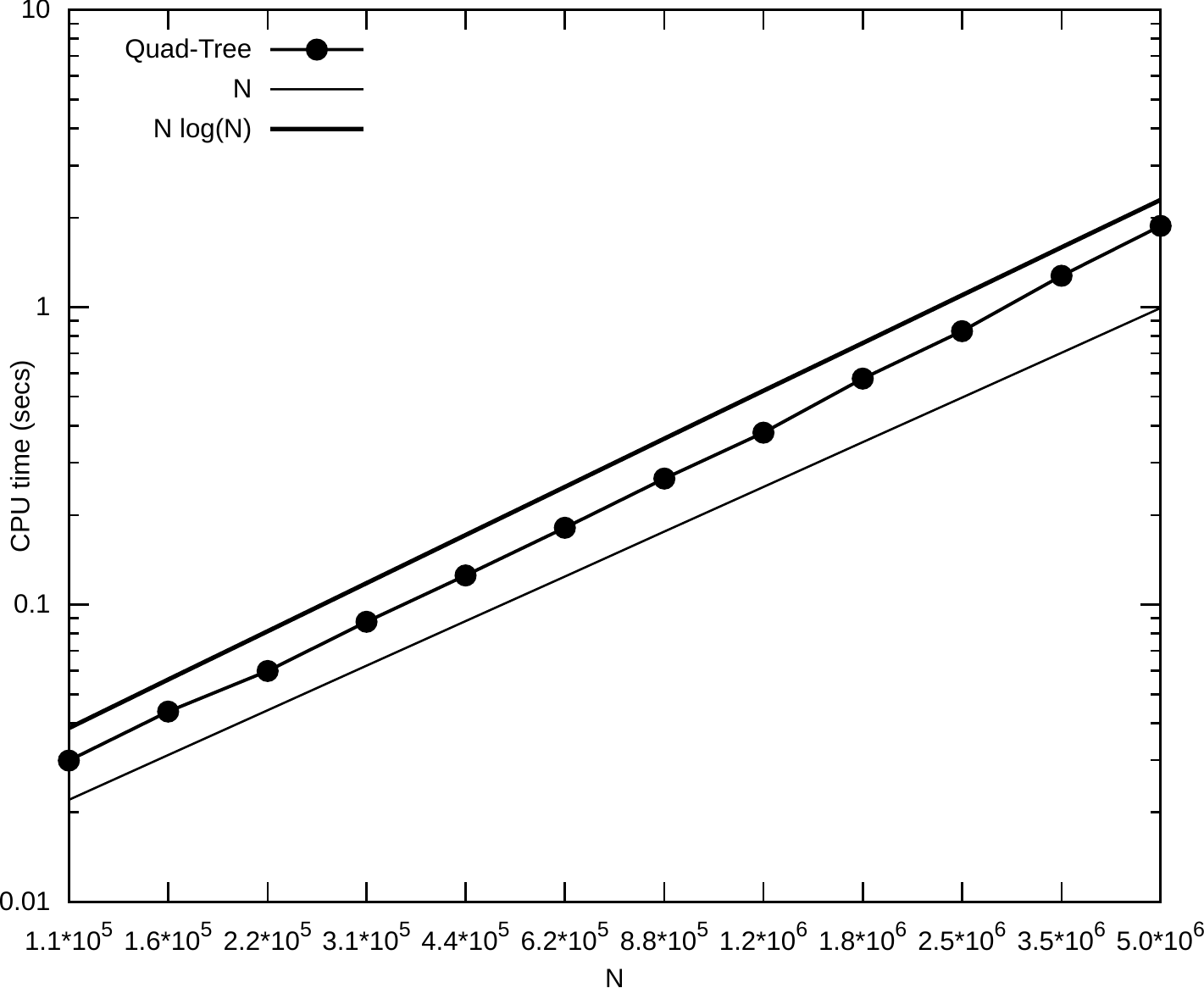}
	\caption{Execution times for the quadtree search on the whole mesh, for $10^5\le N\le 5\cdot 10^6$.}\label{fig:qt_times}
	\end{figure}

Next, we report in Fig. \ref{fig:qt_times} the search times for $N$ random query points, with mesh size $N$ ranging from about $10^5$ to about $5\cdot 10^6$, compared to both $N$ and $N\log N$ orders. As we are performing $N$ searches, each one of expected complexity $\mathcal{O}(\log N)$, the expected order for the total CPU time is $\mathcal{O}(N\log N)$. 
In practice, while this asymptotic behaviour is confirmed, the intermediate scenario can be somewhat less predictable. For example, Fig. \ref{fig:qt_times} shows an almost linear behavior even for a relatively large number of nodes (about $5\cdot 10^5$). This occurs because, under the subdivision rules described above, the resulting, unbalanced quadtree structure reaches its maximum depth only in few regions, compared to the whole mesh (see again Fig. \ref{fig:quadtree}).

\subsection{Barycentric walk} 
\begin{figure}
	\centering
	\includegraphics[height=5cm]{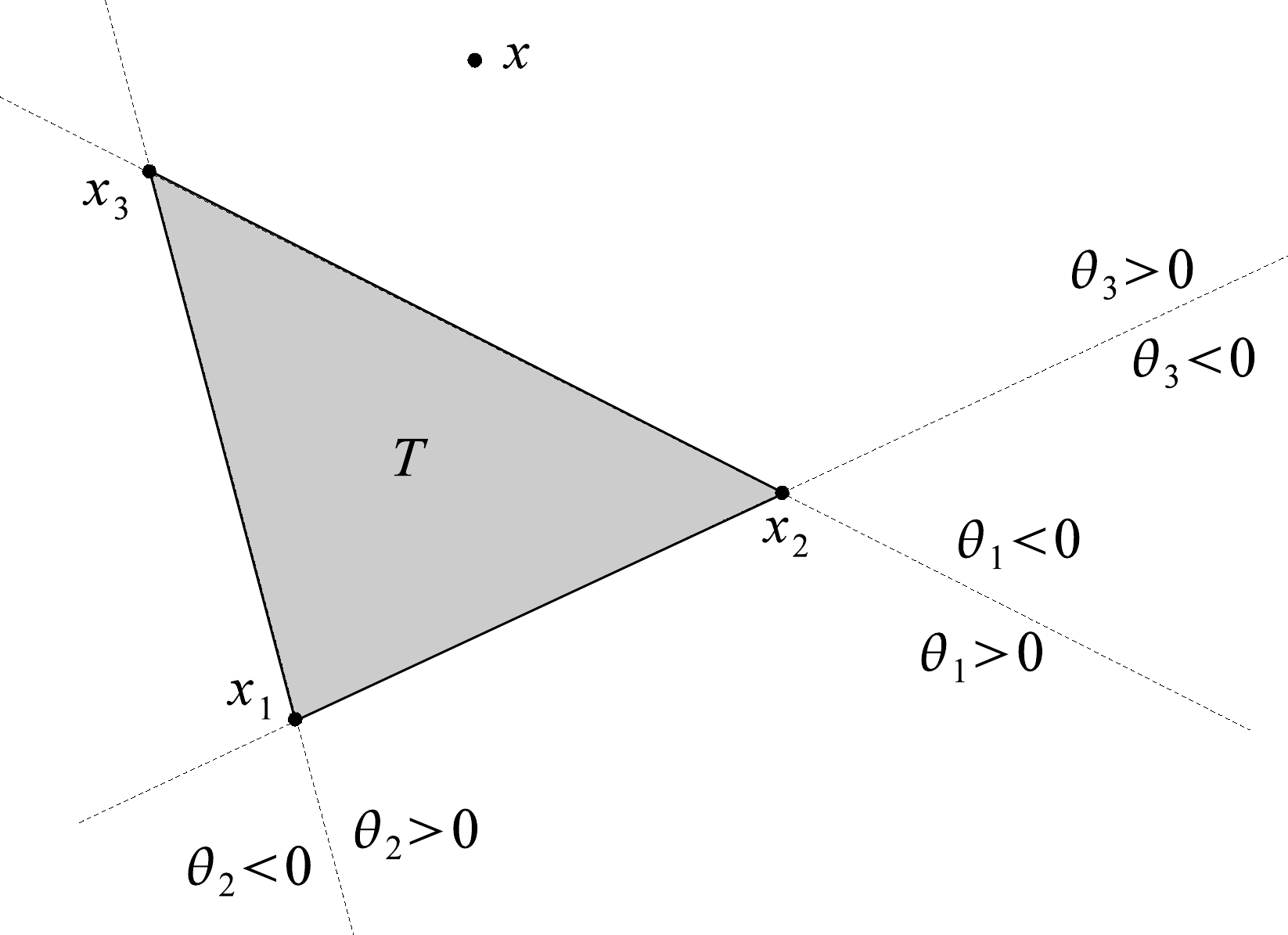}
	\caption{Change of element on the basis of barycentric coordinates.}\label{fig:barycentric}
	\end{figure}

Among the various algorithms which locate a point by stepping along the elements of the triangulation, we review here the so-called {\it barycentric walk}, which is probably the simplest one -- complexity issues are in all cases similar for all the algorithms of this class (see \cite{DPT01} for an extensive review). In this algorithm, in order to locate the point $x$, we start from a given element of the triangulation and change element on the basis of the barycentric coordinates of $x$ with respect to the current element, as shown in Fig. \ref{fig:barycentric}. Given the nodes $x_1$, $x_2$ and $x_3$ of the element $T$ (with $x_i=(\xi_i,\eta_i)$), we write $x=(\xi,\eta)$ by means of the barycentric coordinates $\theta_1$, $\theta_2$, $\theta_3$ as
\[
x = \theta_1 x_1 + \theta_2 x_2 + \theta_3 x_3,
\]
with the $\theta_i$ given by
\begin{equation}\label{eq:barycoords}
\begin{cases}
\theta_1=\frac{(\eta_2-\eta_3)(\xi-\xi_3)+(\xi_3-\xi_2)(\eta-\eta_3)}{(\eta_2-\eta_3)(\xi_1-\xi_3)+(\xi_3-\xi_2)(\eta_1-\eta_3)} \\
\theta_2=\frac{(\eta_3-\eta_1)(\xi-\xi_3)+(\xi_1-\xi_3)(\eta-\eta_3)}{(\eta_2-\eta_3)(\xi_1-\xi_3)+(\xi_3-\xi_2)(\eta_1-\eta_3)} \\
\theta_3=1-\theta_1-\theta_2
\end{cases}
\end{equation}
and we repeat the following steps:
\begin{enumerate}
\item[$a)$]
if all the barycentric coordinates are nonnegative, then $x\in T$ and the point location is complete;
\item[$b)$]
if there exists (at least) one negative coordinate, we look for the node associated with the negative coordinate of largest magnitude, then change element passing to the triangle adjacent to the opposite side, and repeat the computation of the barycentric coordinates on the new element.
\end{enumerate}

For example, in the case shown in Fig. \ref{fig:barycentric}, the only negative coordinate is $\theta_1$, so that we change element from $T$ to the triangle having in common with $T$ the $x_2x_3$ side, see Section \ref{sec:new} for implementation details.

\begin{remark}
The point location ends as soon as we are in an element where all the barycentric coordinates are nonnegative. The finite element-type interpolation on this last element can be immediately computed in terms of these parameters, which are invariant with respect to affine transformations of the reference element. For example, in the $\mathbb{P}_1$ case, we have
\[
I[V](x) = \theta_1 v_1 + \theta_2 v_2 + \theta_3 v_3.
\]
Note also that, whichever algorithm is used for locating the feet of characteristics, the interpolation phase requires to compute the barycentric coordinates in order to interpolate. Therefore, in what follows, the comparison among the various recipes will always include this  computation.
\end{remark}

Concerning the complexity, each change of element has a constant cost, and we can assume that, on a regular Delaunay mesh, the number of elements visited during a walk is asymptotically proportional to its length (this is not necessarily true on graded or anisotropic meshes). Therefore, if we want to locate a query point $Q$ starting the search from a point $P$ (i.e., from a triangle containing this point), the number of walk steps is
\begin{equation}\label{eq:n_steps}
\mathcal{O}\left(\|Q-P\|\sqrt{N}\right),
\end{equation}
where $\sqrt{N}$ is inversely proportional to the space scale of the triangulation, while the location of $Q$ has a complexity of the order of
\begin{equation}\label{eq:bw_times}\mathcal{O}\left(C_1^B+C_2^B\|Q-P\|\sqrt{N}\right),
\end{equation}
in which the constant term accounts for operations which cannot be avoided even in case of a very small distance $\|Q-P\|$: at least one computation of the barycentric coordinates, and, possibly, some change of element. In particular, we observe that if $P$ is fixed, the complexity will be heavier than the quadtree search ($\sqrt{N}$ versus $\log N$). 
We validate our complexity analysis with the following numerical test.

\begin{figure}
	\centering
	\includegraphics[height=6cm]{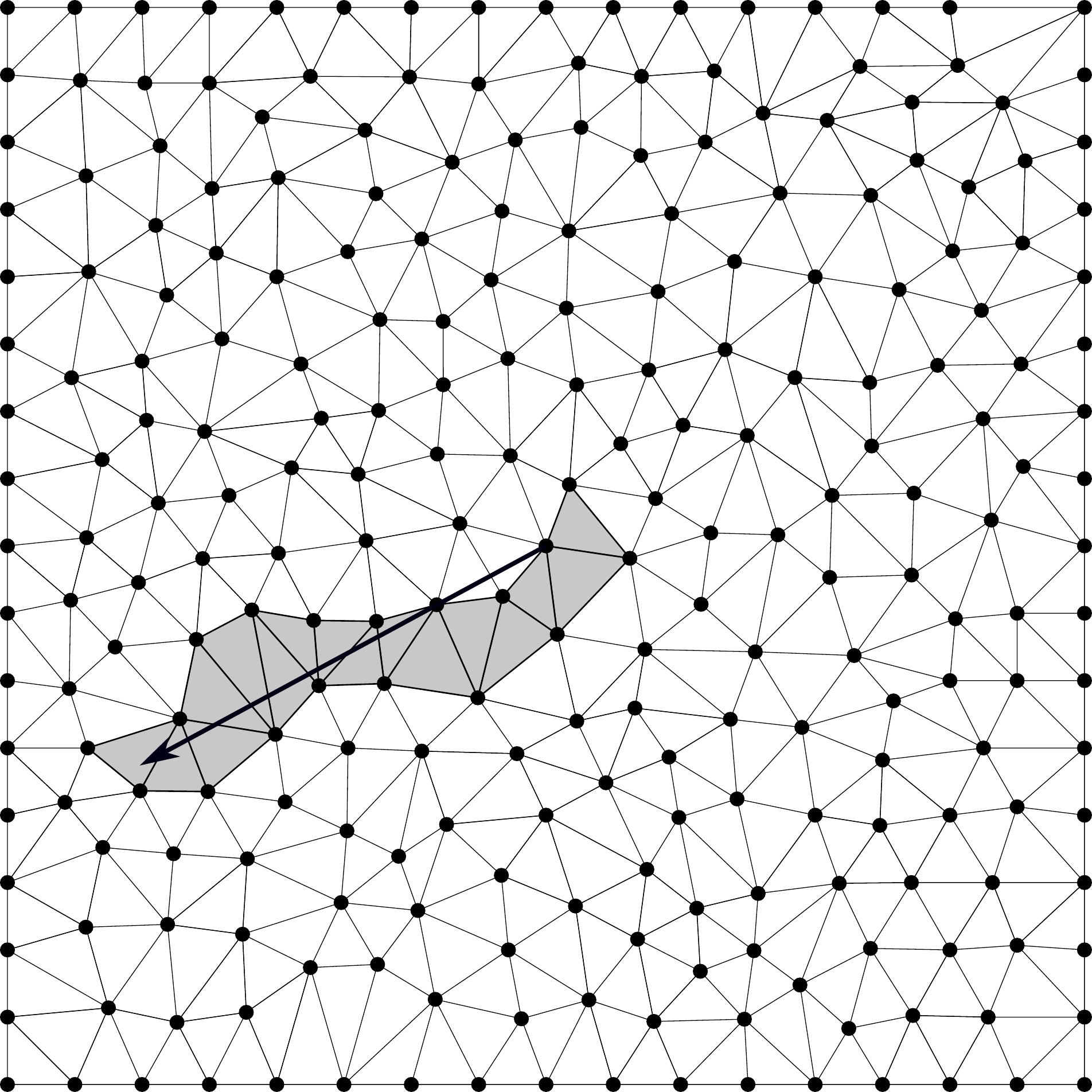}
	\caption{Locating a point via barycentric walk (gray triangles) on a Delaunay triangulation.}\label{fig:bw_steps}
	\end{figure}

\paragraph{Barycentric walk: numerical example.} We first show, on a rough mesh of $N=250$ nodes, the typical barycentric walk for the location of a query point starting from a mesh node, see Fig. \ref{fig:bw_steps}. The initial triangle is randomly chosen from those containing the starting node. 
Note that the barycentric walk is forced, by construction, to perform a large number of steps around those nodes that lie on (or are close to) the line connecting the start and end points. This effect might locally increase the number of steps of the walk, although, as we will soon show, the average number retains a linear dependence on the distance.
	
\begin{figure}
	\centering
	\includegraphics[width=0.47\textwidth]{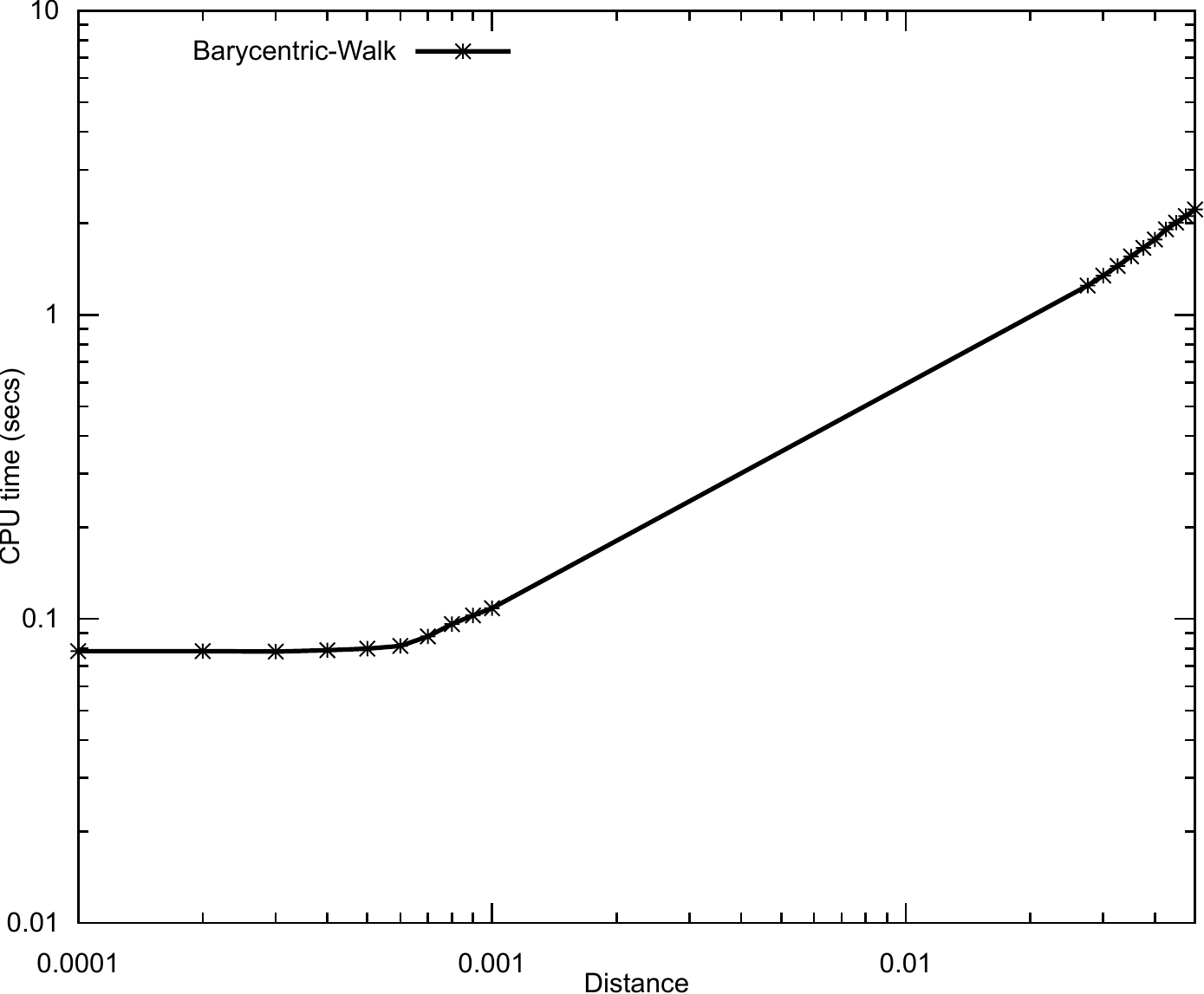} \quad
	\includegraphics[width=0.47\textwidth]{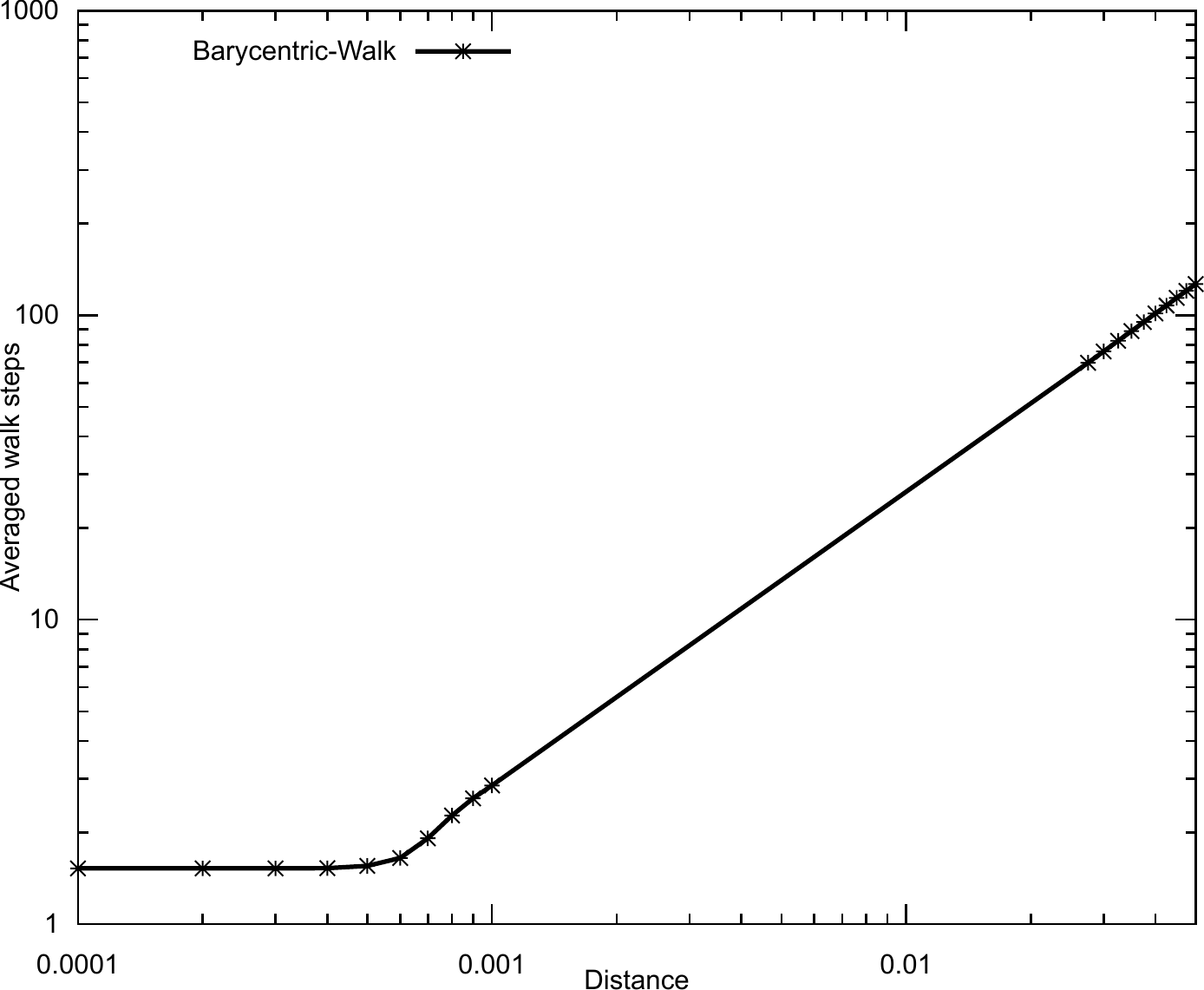}
	\caption{Execution times (left) and averaged walk steps (right) for the barycentric walk search on the whole mesh, versus distance between start and end points, for $N=1.5\cdot 10^6$.}\label{fig:bw_distance}
	\end{figure}

Now, we provide an experimental assessment of the complexity of the algorithm in terms of the distance between the start point and the query point. We choose a fixed mesh of $N=1.5\cdot {10}^6$ nodes, and we compute the total CPU time to locate, for each node, a corresponding query point at given distance. The results are reported in Fig. \ref{fig:bw_distance}. For small distances we clearly observe a plateau in the search times, corresponding to the constant term $C_1^B$ in \eqref{eq:bw_times}, while the behaviour is linear, as expected, when the distance increases (the measure distances have been concentrated in the extreme regions to better catch this behaviour).

\begin{remark}\label{rem:CPU_bw}
In this example, the space scale of the triangulation is estimated by $1/\sqrt{N}=8\cdot 10^{-4}$, while we observe that CPU time begins to grow already at a smaller distance. This reflects the fact that the actual intersection between the trajectory and each element amounts in general to a fraction of the space scale, as clearly shown by Fig. \ref{fig:bw_steps}. The behavior of the averaged walk steps is similar, in particular we observe a value of about $1.5$ for the plateau. Here, the random choice of the initial triangle of the walk at each mesh node implies on average some change of element even at a very small distance.
\end{remark}

\subsection{Quadtree and barycentric walk complexity in space} In this section, we briefly compare the quadtree and the barycentric walk in terms of space complexity, i.e., of memory occupation. Recall that we have heuristically assumed that the number of elements is $\mathcal O(N)$, and that the average depth of the quadtree is $\mathcal O(\log N)$. Both point location approaches use the mesh information, namely the list of point coordinates of each grid node and the connectivity, in the form of a list of triplets of vertex indices, ordered as they appear in the node list. 
In addition, the barycentric walk requires, for each triangle, the list of triangle neighbors to move across the elements. This list consists in triplets of triangle indices ordered as they appear in the connectivity list. On the other hand, the quadtree structure is more complicate. Starting from the root, each node of the tree must record the four coordinates (for the left/bottom and top/right vertices) of its quad, and four pointers to its children, while the leaves contain the indices of mesh nodes and triangles intersected by their quads. According to the rules discussed in Sect. \ref{QT-rules}, the construction stops if a quad contains at most one mesh node, regardless the number of incident triangles, and at most $q$ triangles if it contains no mesh nodes. This implies that the total number of vertex indices in the leaves is about $N$ (some duplicates can be found if a mesh node stands on the side or is exactly one vertex of a quad), whereas the total number of triangle indices is much greater than the number of triangles, since a triangle typically overlaps with several quads. Note that, for point location, only the triangle indices are needed. Hence, in the following computation, we drop the list of vertex indices after the quadtree construction.

In order to estimate the order of memory occupation for the quadtree, we note that, starting from the root, each successive level has four times the number of quads of the previous. The total memory occupation is proportional to the total number of quads, i.e.,
\[
\sum_{k=0}^{\log N} 4^k = \mathcal O(N),
\]
as it can be easily seen, for example, via comparison with an integral. On the other hand, the barycentric walk requires to store the list of neighboring elements for each triangle, which results again in a linear memory occupation.

The following numerical test validates the expected $\mathcal{O}(N)$ space complexity in terms of storage for the corresponding data structures.

\paragraph{Storage for quadtree and barycentric walk: numerical example.} 

\begin{figure}
	\centering
	\includegraphics[width=0.47\textwidth]{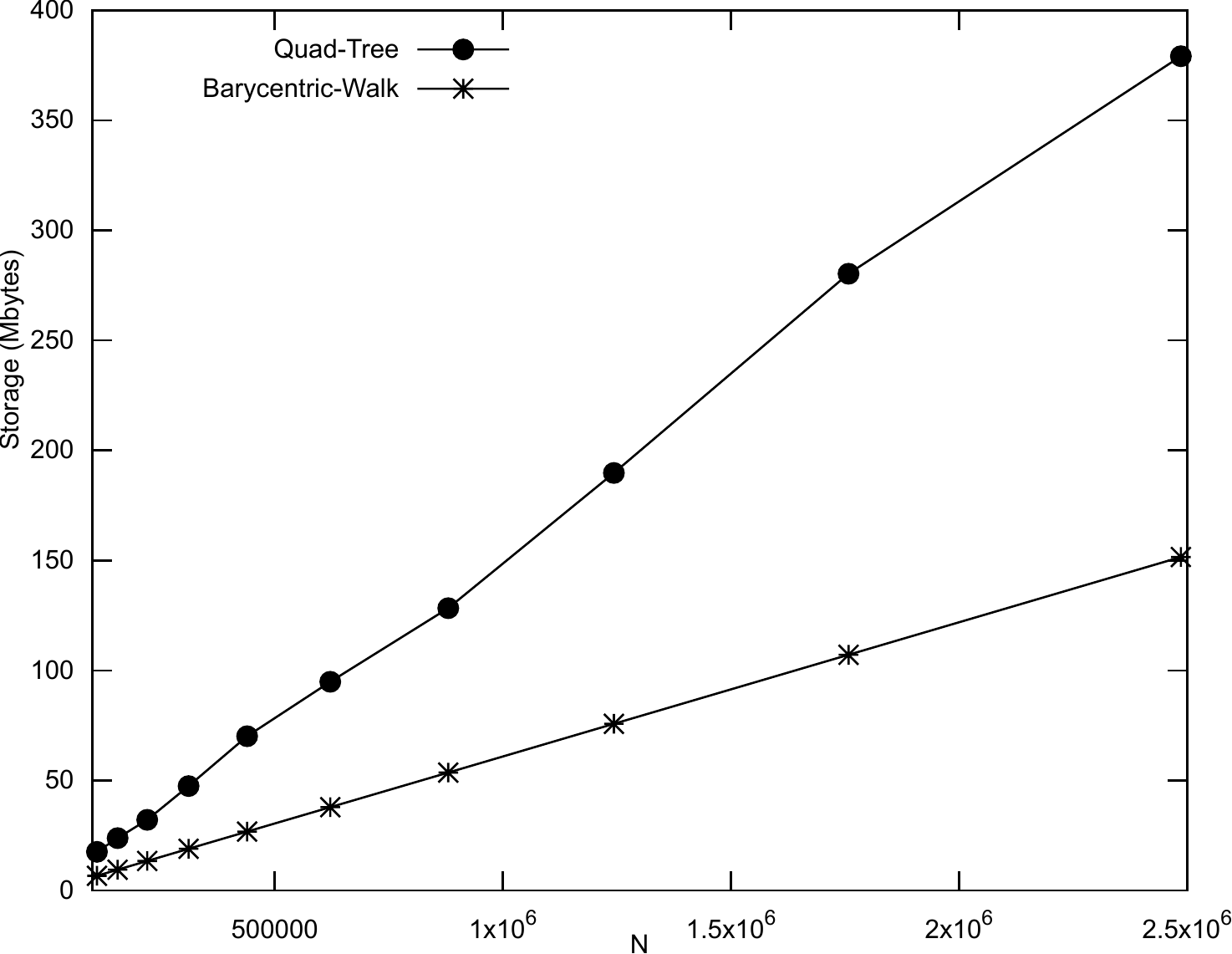} 
	\caption{Space complexity of quadtree and barycentric walk in terms of storage (Mbytes), for $10^5\lesssim N\lesssim 2.5\cdot 10^6$.}\label{fig:qt-bw_storage}
	\end{figure}

We report in Fig. \ref{fig:qt-bw_storage}, for meshes of size $10^5\lesssim N\lesssim 2.5\cdot 10^6$, the storage in Mbytes corresponding to the two data structures, including for both the load due to the mesh (vertices plus triangles). The experiment confirms the $\mathcal O(N)$ space complexity for both approaches. Nevertheless, we found that the number of tree nodes is about $2N$, while the number of triangle indices in the leaves is about $9N$, and this results in a gain factor about $2.5$ for the barycentric walk. For completeness, we remark the code has run on a $64$bit architecture, in which the storage for integers, doubles and pointers amounts, respectively, to $4$, $8$ and $8$ bytes each.  

\section{Element search by barycentric walk}\label{sec:new}

While the quadtree search starts necessarily from the same root at each execution, the barycentric walk, which has in the general case a higher complexity, might nevertheless be improved by using a ``good'' choice of the starting element, and keeping memory of such choice. This is what will be pursued in this section, and is the key to obtain an $\mathcal{O}(1)$ complexity for the location of a point on the grid, when tracking characteristics via a SL scheme for advection equations. Such a reduction of complexity makes it convenient to replace the quadtree algorithm, which is memory-consuming and complex to code, with an easier and smarter procedure.

To fix ideas, we consider a regular Delaunay triangulation with $N$ nodes and space scale $\Delta x\propto 1/\sqrt{N}$, and we use Euler tracking of characteristics \eqref{eq:eulero}, that we recall here for the reader's convenience:
$$X^\Delta(x_i,t_{n+1};t_n) = x_i-\Delta t \> f(x_i,t_{n+1}).$$
In this setting, the complexity \eqref{eq:bw_times} for the barycentric walk algorithm reads
$$
\mathcal{O}\left(C_1^B+C_2^B\frac{\|X^\Delta(x_i,t_{n+1};t_n)-\overline X_i\|}{\Delta x}\right),
$$
where $X^\Delta$ is the query point and $\overline X_i$ is a starting point (or the corresponding starting element) related to the node $x_i$ from which the characteristic originates. Now, we analyze some choice of $\overline X_i$, suitable to obtain a point location with a complexity independent of the grid size. As a first choice, we consider the strategy of tracking all the characteristic. This technique has already been applied to SL schemes, as discussed in the introductory overview.

\begin{enumerate}
	
\item[$a)$] $\overline X_i=x_i$. This strategy is used, for example, in \cite{RBS06,BDR13,B20}, coupled with a substepping along the characteristic. With this choice, 
$$\|X^\Delta(x_i,t_{n+1};t_n)-x_i\|=\Delta t\|f(x_i,t_{n+1})\|,$$
so that the number of steps is of the order of the local Courant number, and a single element search has therefore a complexity of
$$
\mathcal{O}\left(C_1^B+C_2^B\frac{\|f(x_i,t_{n+1})\|\Delta t}{\Delta x}\right),
$$
that is, asymptotically constant under linear $\Delta t/\Delta x$ relationship. Note that, if one works at large Courant numbers in order to increase efficiency of the scheme, the element search becomes in turn more complex. Moreover, complexity is no longer asymptotically constant under nonlinear refinements in which $\Delta x=o(\Delta t)$.
	
\item[$b)$] $\overline X_i=X^\Delta\left(x_i,t_n;t_{n-1}\right)$. Note that this point is the foot of the characteristic at the previous time step (alternatively, the element containing this point), and has already been computed. In this case,
\begin{eqnarray*}
\left\|X^\Delta(x_i,t_{n+1};t_n)-X^\Delta(x_i,t_n;t_{n-1})\right\| & = & \Delta t \|f(x_i,t_{n+1})-f(x_i,t_n)\| \\
& \le & L_t \Delta t^2,
\end{eqnarray*}
and the location of the foot of characteristics has therefore a complexity of
$$
\mathcal{O}\left(C_1^B+C_2^B L_t \frac{\Delta t^2}{\Delta x}\right),
$$
in which, since $L_t$ is a global Lipschitz constant, we are bounding the computational cost from above.
In this case, the complexity is asymptotically constant provided $\Delta t=\mathcal{O}\left(\Delta x^{1/2}\right)$. On the other hand, under a linear refinement, it tends to coincide with the complexity of a single computation of the barycentric coordinates: in other terms, the event of a change of element becomes more and more unlikely. In particular, regions of the domain in which the advecting vector field has slow changes (or tends towards a regime state) require only minor adjustments from one time step to the next. In the limit case of an advection term constant in time, no change of element is necessary.\\
Note that this initialization is constructed independently for each node, and hence the location of the points $X^\Delta(x_i,t_{n+1};t_n)$ can be performed in parallel w.r.t. $i$. Since it requires the same sequence of operations for each node (except for a possibly different length of the walk), the resulting algorithm might be particularly convenient on a SIMD architecture.

\item[$c)$] $\overline X_i=X^\Delta\left(x_k,t_{n+1};t_n\right)$, namely the foot of the characteristic at the same time, but at a node $x_k$ adjacent to $x_i$. With this choice,
\begin{eqnarray*}
\|X^\Delta(x_i,t_{n+1};t_n)-X^\Delta(x_i,t_{n+1};t_n)\| & = &  \|x_i-x_k-\Delta t(f(x_i,t_{n+1})-f(x_k,t_{n+1}))\| \\
& \le & (1+L_x \Delta t)\Delta x.
\end{eqnarray*}
Taking into account the fixed complexity terms, the point location has therefore a cost of
$$
\mathcal{O}\left((C_1^B+C_2^B)+C_2^B L_x \Delta t \right).
$$
Again, we obtain an asymptotically constant complexity, but it appears to have a less critical dependence (if any dependence at all) on the $\Delta t/\Delta x$ relationship, and in particular to be applicable when $\Delta x=o(\Delta t)$. However, opposite to what happens in the previous case $(b)$, even with stationary advection terms, we expect that a change of element is needed in general, and this causes an increase of the constant term.

This strategy clearly requires that the nodes are put in a sequence where all nodes but the first one have a neighbour for which the final element has already been computed. In practice, this may be accomplished by constructing a spanning tree of the grid, once and for all after the grid construction. Note that, with respect to the previous strategy, this technique is more complicate to set in parallel form. Parallelization should be performed on successive levels of the spanning tree, by computing in parallel all the nodes having parents at the previous level, and its efficiency is clearly related to the depth of the spanning tree, and ultimately to the mesh size.

Note also that, in SL schemes (see, e.g., the discussion of this point in \cite{FM20}), it is usually required for stability reasons that characteristics passing through neighbouring nodes do not cross. In practice, $X^\Delta(x_i,t_{n+1};t_n)$ and $X^\Delta(x_k,t_{n+1};t_n)$ must always have a positive distance, so that, using \eqref{eq:eulero} and the reverse triangular inequality,
\begin{eqnarray*}
\|X^\Delta(x_i,t_{n+1};t_n) - X^\Delta(x_k,t_{n+1};t_n)\| & \ge & \|x_i-x_k\| - \Delta t \|f(x_i,t_{n+1})-f(x_k,t_{n+1})\| \\
& \ge & (1-\Delta t \>L_x)\|x_i-x_k\| > 0.
\end{eqnarray*}
This leads to the well-known condition
\[
L_x\Delta t<1,
\]
and, as a consequence,
\begin{equation}\label{eq:uniform-bound}
\|X^\Delta(x_i,t_{n+1};t_n)-X^\Delta(x_k,t_{n+1};t_n)\| \le 2\Delta x,
\end{equation}
which also implies a uniform bound w.r.t. $\Delta x$ on the complexity.
\end{enumerate}
\begin{remark}
For an actual implementation, all the three walk strategies presented above require some additional data structures with respect to the standard barycentric walk. More precisely, we need a list of $N$ integers for storing the indices of all the triangles, one for each vertex of the mesh, from which to start the barycentric walks. Furthermore, strategy $(c)$ also  requires the spanning tree for ordering the grid nodes. This results in another list of $N$ integers, storing for each node the index of its parent in the spanning tree. In the left plot of Fig. \ref{fig:bwabc_storage} we report, for meshes of size  $10^5\lesssim N\lesssim 2.5\cdot 10^6$, the storage in Mbytes corresponding to the different strategies (clearly the same for strategies $(a)$ and $(b)$), including the the mesh data and the neighbor list for the standard barycentric walk. Finally, the right plot reports the storage improvement, showing that the walk strategies require less than half the memory required by the quadtree location.  

	\begin{figure}
	\centering
	\includegraphics[width=0.47\textwidth]{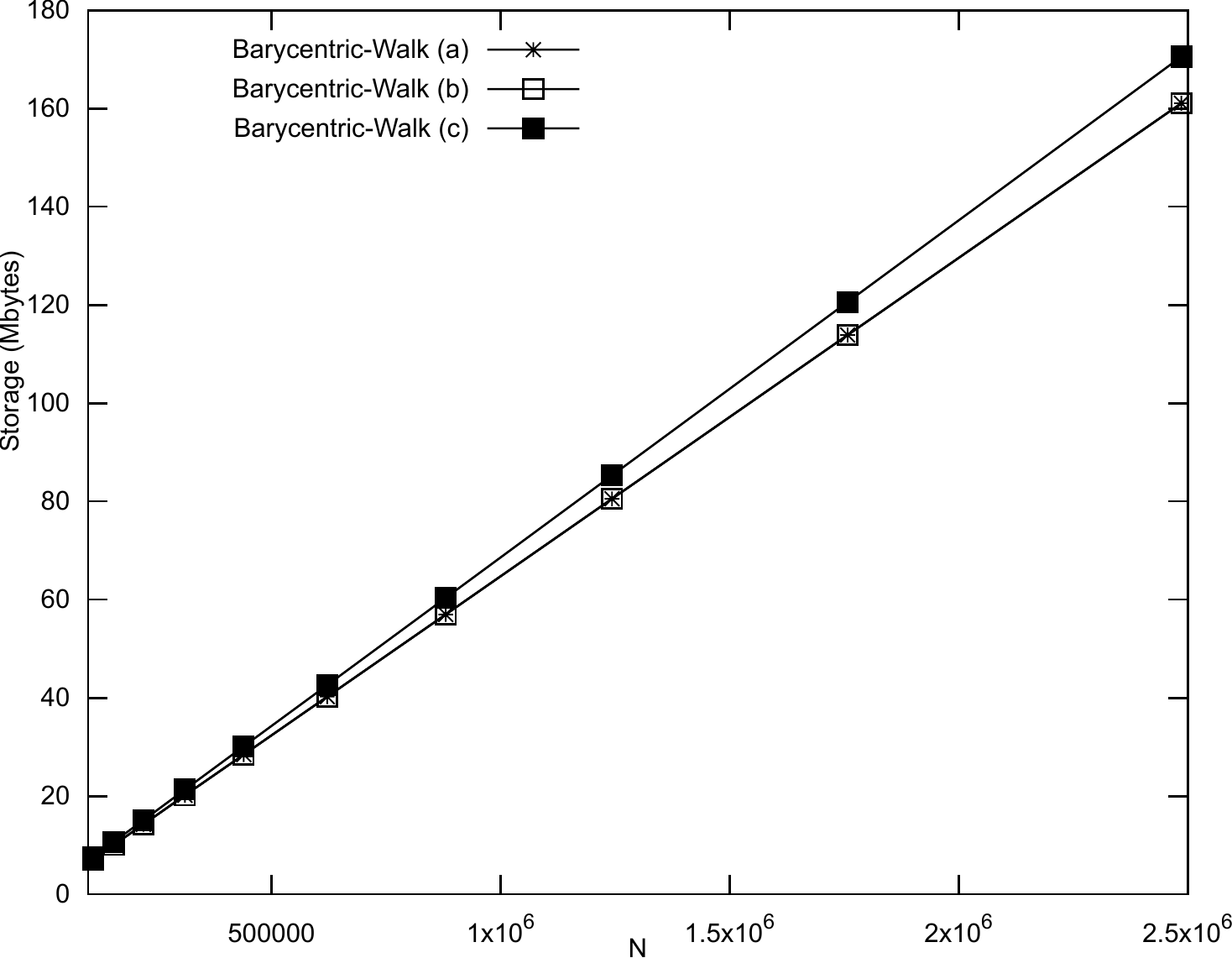} \quad
	\includegraphics[width=0.47\textwidth]{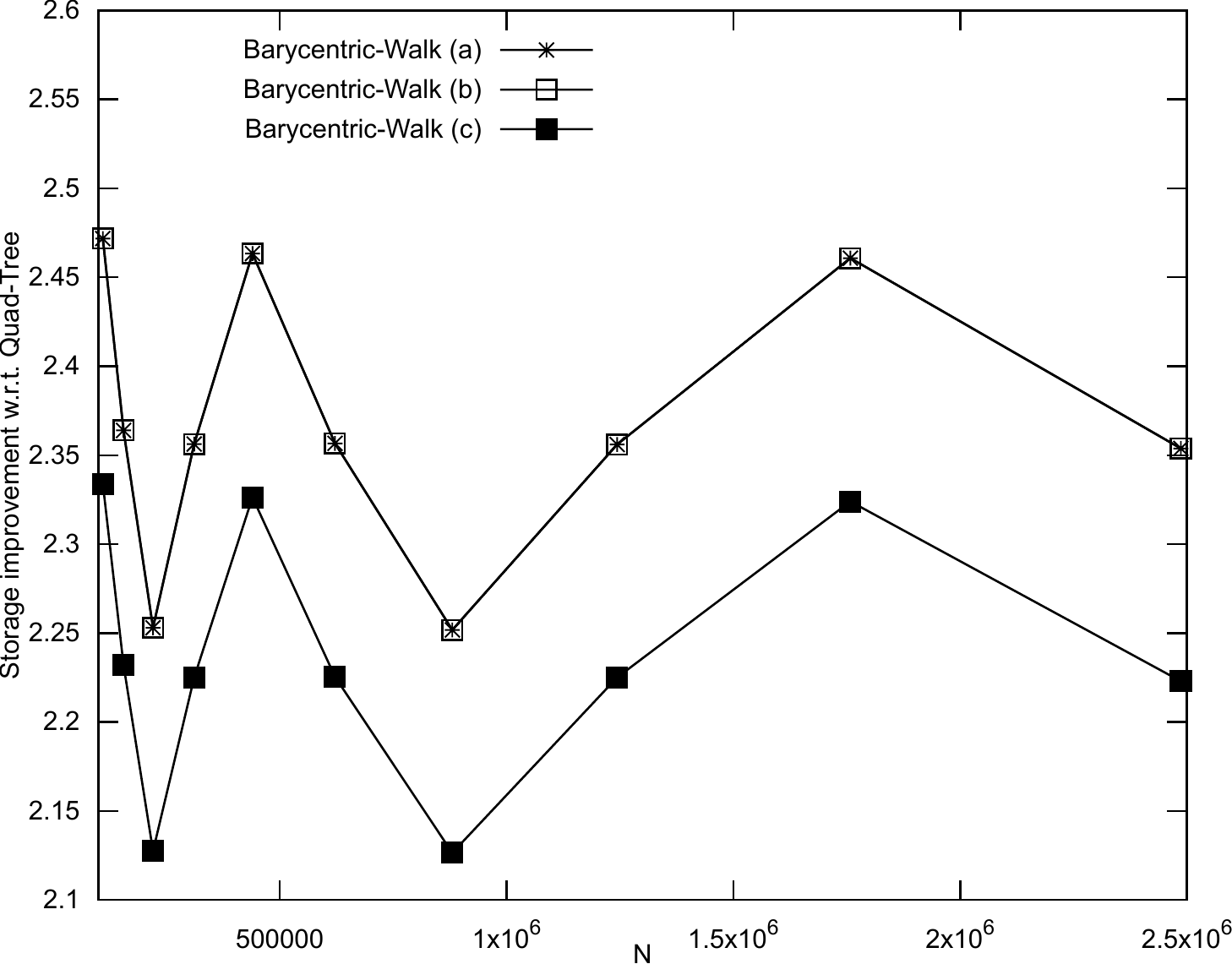}
		\caption{Storage (Mbytes) for the walk strategies (left) and improvement w.r.t. storage for the quadtree, for $10^5\lesssim N\lesssim 2.5\cdot 10^6$.}\label{fig:bwabc_storage}
	\end{figure}
\end{remark}

To conclude this section, we provide some implementation details and a pseudo-code for our barycentric walk algorithm. We consider, as follows, a suitable data structure containing all the relevant information for the triangulation. We recall that $\mathcal{V}=\{x_i\}$, for $i=1,\dots,N$, is the list of point coordinates of each node. We denote by $\mathcal{T}$ the list of triangles, namely the list of vertex indices $\mathcal T_j=(i_1,i_2,i_3)\in\{1,\dots,N\}^3$ defining the triangle with vertices $x_{i_1},x_{i_2},x_{i_3}$, for $j=1,\dots,N_t$, where $N_t$ is the total number of triangles. Moverover, we denote by $\mathcal N$ the list of triangle neighbors, namely the list of triangle indices $\mathcal N_j=(j_1,j_2,j_3)\in\{1,\dots,N_t\}^3$ corresponding to the three neighbors $\mathcal T_{j_1},\mathcal T_{j_2},\mathcal T_{j_3}$ of $\mathcal T_{j}$, for $j=1,\dots,N_t$. We adopt the standard convention for which the index $j_k$ corresponds to the neighboring triangle of $\mathcal T_j$ sharing the edge opposite to the vertex with index $i_k$, for $k=1,2,3$. 
Furthermore, we denote by $\mathcal{T}^0$ a list of triangle indices associated to the nodes, so that $\mathcal{T}^0_i\in\{1,\dots,N_t\}$, for $i=1,\dots,N$, identifies the initial triangle for the barycentric walk which tracks the characteristic originating from the node $x_i$. We always initialize $\mathcal{T}^0$ assigning to each node $x_i$ a random triangle among those having $x_i$ as a vertex.
We remark that, to implement the walk strategy $(c)$, we need a root node for the spanning tree of the grid. For simplicity, we assume that $\mathcal{V}$ is already ordered according to the spanning tree, so that $x_1$ is the root node, followed by its first neighboring nodes, and so on recursively. Then we denote by $\mathcal{P}$ the list of indices of parent nodes, so that $\mathcal P_i\in\{1,\dots,N\}$, for $i=1,\dots,N$, identifies the parent node of $x_i$. In particular, the root node is the only one satisfying $\mathcal P_1=1$. Finally, given the dynamics $f$, we can build the list $\mathcal{Q}^n=\{q_i^n\}_{i=1,\dots,N}$ of query points for the barycentric walk at time $t_n=n\Delta t$, tracking the characteristics with a suitable solver for ordinary differential equations (e.g., $q^n_i=x_i-\Delta t\,f(x_i,t_n)$ for the Euler scheme). The procedure is implemented as a pseudo-code in Algorithm \ref{ALG1}.

\begin{algorithm}[!t]
\small
\begin{algorithmic}[1]
\State Given $\mathcal V$, $\mathcal T$, $\mathcal N$, $\mathcal T^0$, $\mathcal P$, $\mathcal Q^n$, and a strategy $s\in\{a,b,c\}$.
\State {\bf function} $(\mathcal T^f,\mathcal B)=\mbox{PointLocationBW}(\mathcal Q^n,s)$
\For{$i=1:N$}
 \If{$s\neq c$}
 \State Set $j\leftarrow \mathcal{T}^0_i$
 \Else
 \State Set $j\leftarrow \mathcal{T}^0_{\mathcal P_i}$
 \EndIf
 \State Set BW $\leftarrow$ true
 \While{BW}
\State Set $(i_1,i_2,i_3)\leftarrow \mathcal{T}_{j}$
\State Compute $(\theta_1,\theta_2,\theta_3)$ for $q^n_i$ w.r.t. $x_{i_1}, x_{i_2}, x_{i_3}$ using \eqref{eq:barycoords} 
\State Set $\theta^*\leftarrow \displaystyle\min_{k=1,2,3}\theta_k$
and $k^*\leftarrow\arg\displaystyle\hskip-7pt\min_{k=1,2,3}\theta_k$
\If{$\theta^*\ge0$}
\State Set BW $\leftarrow$ false
\Else
\State Set $(j_1,j_2,j_3)\leftarrow \mathcal{N}_{j}$
\State Set $j\leftarrow j_{k^*}$ 
\EndIf
\EndWhile 
 \State Set $\mathcal T^f_i\leftarrow j$
 \State Set $\mathcal{B}_i\leftarrow (\theta_1,\theta_2,\theta_3)$
 \If{$s=c$}
 \State Set $\mathcal T^0_i\leftarrow j$
 \EndIf
\EndFor
\If{$s=b$}
 \State Set $\mathcal T^0\leftarrow \mathcal T^f$
 \EndIf
\State \Return $\mathcal T^f$ and $\mathcal B$
\State {\bf end function} 
\end{algorithmic}
\caption{Pseudo-code for the function {\tt PointLocationBW}}\label{ALG1}
\end{algorithm}
We remark that, when the code is meant to run for several time steps in a scheme for advection equations, the walk strategy $s\in\{a,b,c\}$ should be chosen once and for all at the beginning, since both strategies $(b)$ and $(c)$ overwrite the list $\mathcal{T}^0$. In particular, for the very first time step, the walk strategies $(a)$ and $(b)$ coincide. Here, we implemented different switches on $s$ just to enlighten the differences between the walk strategies and simplify the presentation. On the contrary, in the following tests, we implemented three different versions of the algorithm, thus saving several {\tt if} statements ($2N$ per time step) at runtime. 

For the sake of completeness, we also report, in Algorithm \ref{ALG2}, the basic SL pseudo-code for the advection equation \eqref{eq:trasp}, using $\mathbb{P}_1$ interpolation for the solution and the Euler scheme for tracking characteristics. This will be used later to show which percentage of the total computational load is due to the point location.
\begin{algorithm}[!t]
\small
\begin{algorithmic}[1]
\State Given $\Omega\subset\R^2$, $u_0:\Omega\to \R$, $f:\Omega\times\R^+\to\R^2$, $\Delta x, \Delta t, T>0$
\State Choose a walk strategy $s\in\{a,b,c\}$
\State Build a triangulation $(\mathcal V,\mathcal T,\mathcal N,\mathcal T^0, \mathcal P)\leftarrow (\Omega,\Delta x)$
\State Set $V^0_i\leftarrow u_0(x_i)$ for $i=1,\dots,N$
\State Set $n\leftarrow 0$
\While{$n\Delta t\le T$}
\State Set $q^n_i\leftarrow x_i-\Delta t\,f(x_i,t_n)$ for $i=1,\dots,N$
\State Compute $(\mathcal T^f,\mathcal B)\leftarrow \mbox{PointLocationBW}(\mathcal Q^n,s)$ using Algorithm \ref{ALG1} 
\For{$i=1:N$} 
\State Set $(i_1,i_2,i_3)\leftarrow \mathcal T^f_{i}$
\State Set $(\theta_1,\theta_2,\theta_3)\leftarrow \mathcal{B}_i$
\State Set $V^{n+1}_i\leftarrow\theta_1 V^n_{i_1}+\theta_2 V^n_{i_2}+\theta_3 V^n_{i_3}$
\EndFor
\State Set $n\leftarrow n+1$
\EndWhile
\end{algorithmic}
\caption{Pseudo-code for the basic Semi-Lagrangian advection}\label{ALG2}
\end{algorithm}

\section{Numerical examples}\label{sec:numeric}
In this section, we present several numerical tests, showing the performance of the proposed search strategies, as compared with the standard quadtree search and with a direct search on a structured grid. Moreover, we provide a comparison with the built-in Matlab function {\tt pointLocation}. Finally, we present the performance of the basic SL scheme, equipped with the different point location algorithms.

All the algorithms have been implemented from scratch in C++ language, compiled with GCC compiler 7.5.0, and run (in serial for the present work) on a  PC Desktop equipped with an Intel i9-9900K CPU with 16 cores 3.60Ghz, 32Gb RAM, under the OS Ubuntu 18.04.3 LTS. We have also built a simple wrapper to easily employ the library {\tt Triangle} for the generation of quality Delaunay meshes \cite{TRI}. In particular, the library accepts an input constraint $A_{\max}$ for the maximal area of each triangle in the mesh. Then, we set the space scale $\Delta x=\sqrt{2A_{\max}}$, so that triangle areas are proportional to $\frac12 \Delta x^2$.

\begin{figure}
	\centering
	\includegraphics[height=6cm]{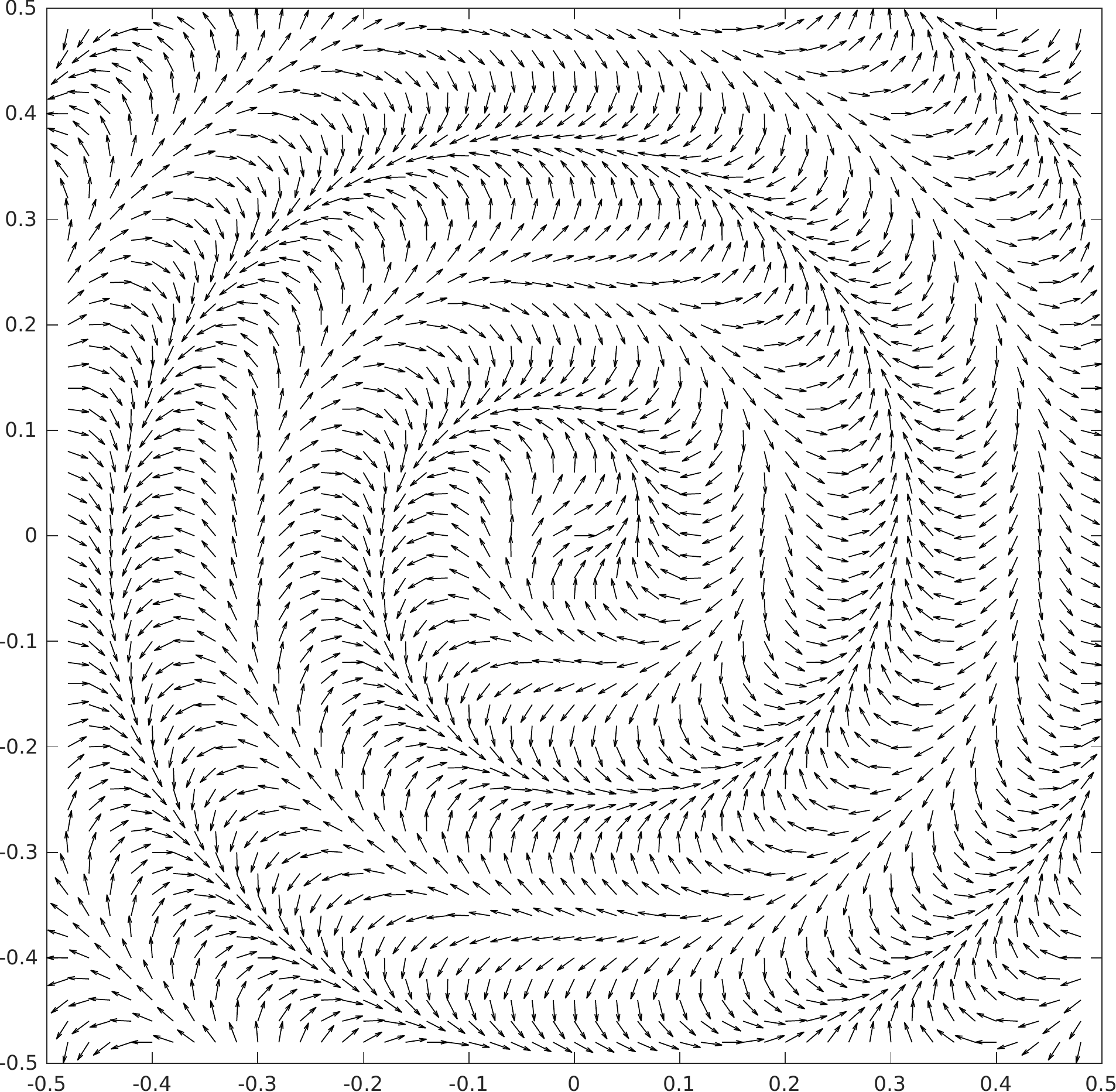}
	\caption{Test dynamics with $C_0=8\pi$, $C_1=2\pi$ at time $t=0$, sampled on a uniform grid with $\Delta t/\Delta x=1$.}\label{fig:dynamics}
	\end{figure}

In all the tests, we consider the following advecting vector field:
\[
f(x,t) = \begin{pmatrix}\cos \left(C_0\|x\|+C_1 t\right) \\ \sin \left(C_0\|x\|+C_1 t\right) \end{pmatrix}\,,
\]
with $x\in\Omega$, $t\in[0,1]$ and $C_0, C_1>0$, namely a rotating vector field with frequencies $C_0$ and $C_1$, respectively in space and time (see Fig. \ref{fig:dynamics}). 
	
 In order to compare the execution times with the structured case, we take the square domain $\Omega=[-1/2,1/2]^2$, and we exclude from the computation all the grid nodes for which the corresponding characteristic, tracked by the Euler scheme \eqref{eq:eulero}, falls out the domain. Moreover, in the construction of the spanning tree for the walk strategy $(c)$, we take the root node as the closest to the center of the domain.

Note that, by definition, we have $\|f\|=1$ everywhere in $\Omega$, while the Lipschitz constants of $f$ are given by $L_x=C_0$ and $L_t=C_1$. This allows to better analyze the complexity of the search strategies in terms of the Courant number 
$\|f\|\Delta t /\Delta x$, which is indeed the same on the whole domain: setting $\Delta t=\alpha\Delta x$, 
with $\alpha\ge 0$, the Courant number is simply given by $\alpha$. Finally, we average in time, dividing by the number of time steps $\ceil*{1/\Delta t}$, both the computational times and the averaged walk steps.

\begin{figure}
	\centering
	\includegraphics[width=0.47\textwidth]{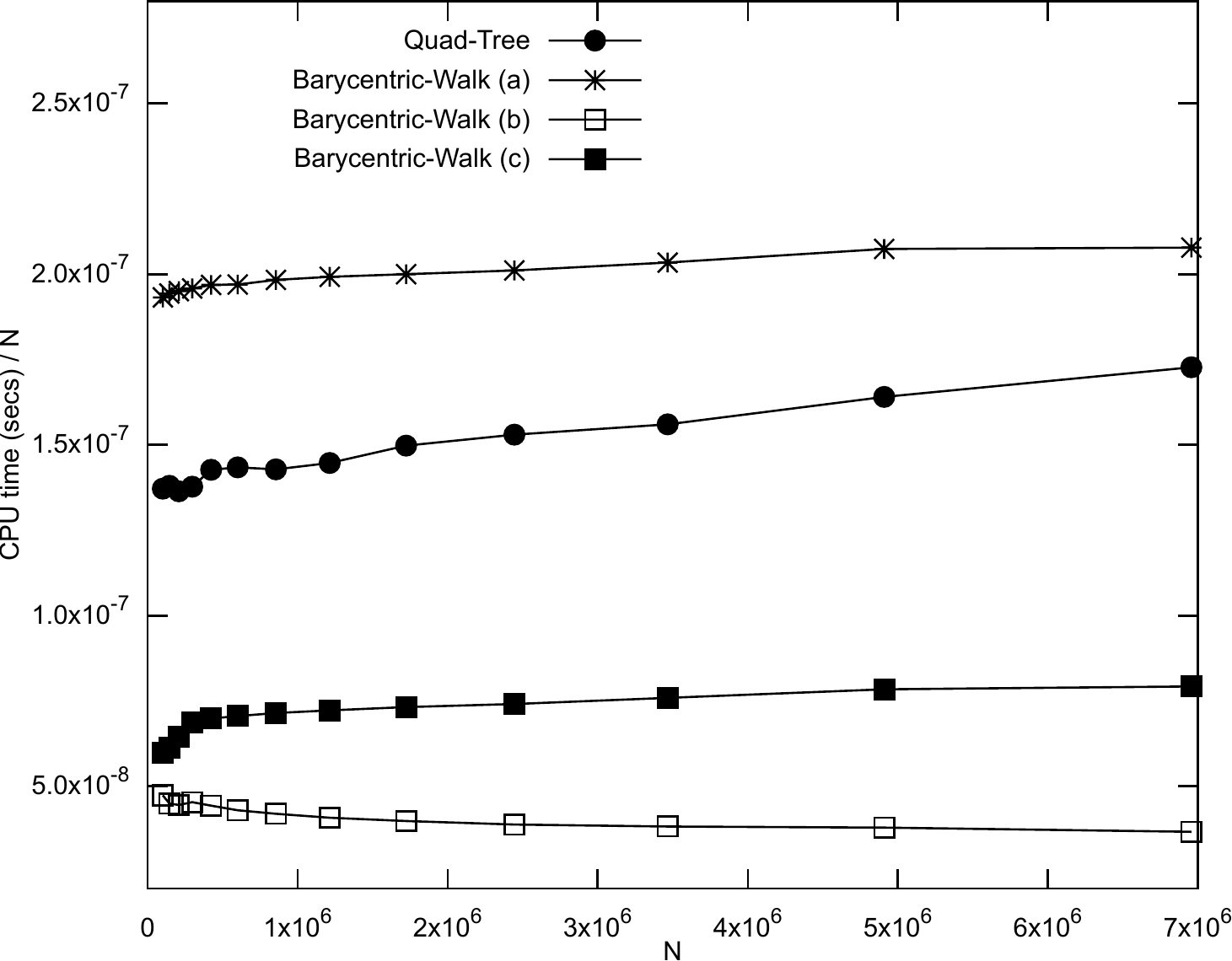} \quad
	\includegraphics[width=0.47\textwidth]{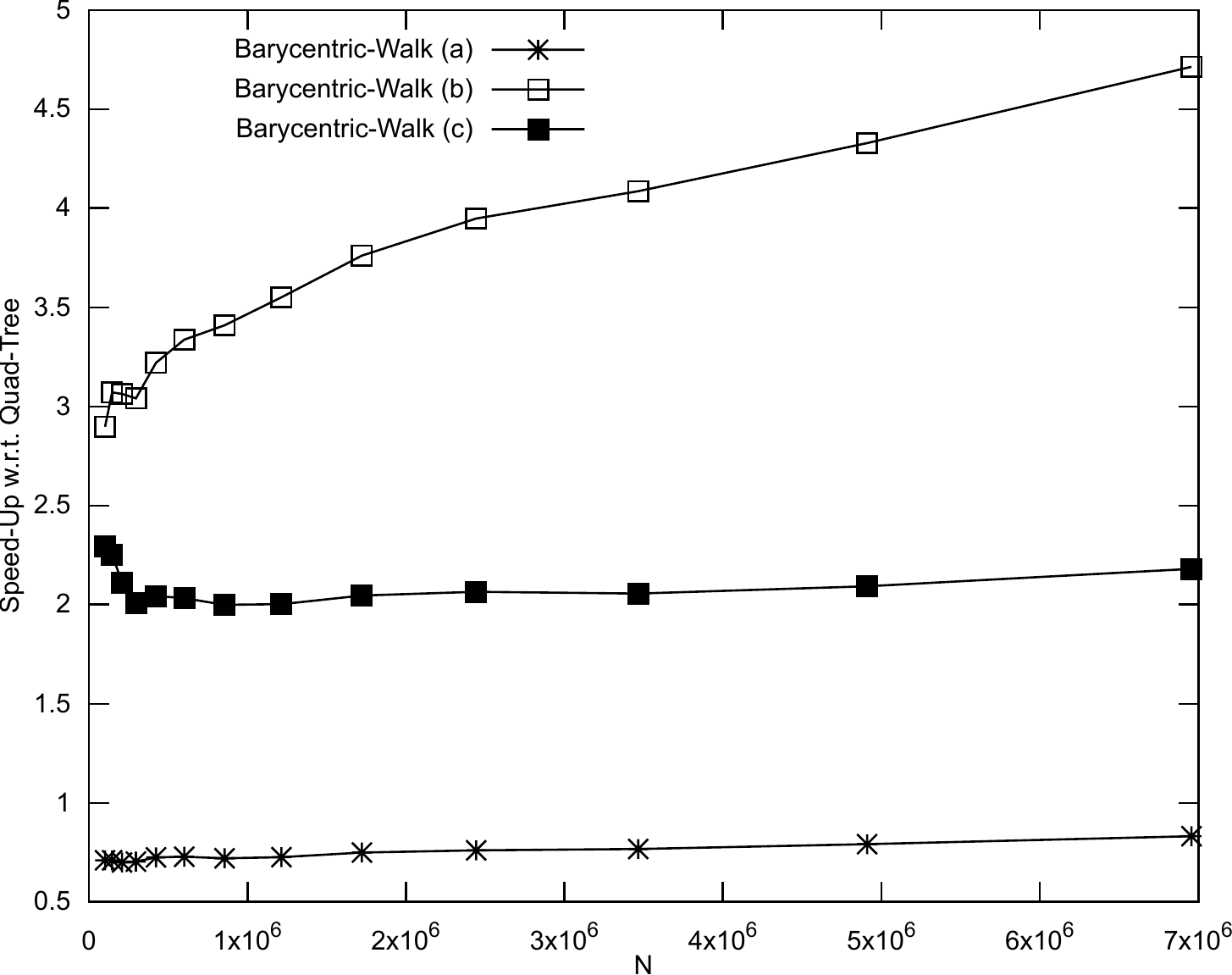}
	\caption{Averaged search time of the various strategies (left) and improvement factor of the barycentric walk strategies versus the quadtree search (right), for $10^5\lesssim N\lesssim 7\cdot 10^6$.}\label{fig:comparison}
	\end{figure}

\paragraph{Comparison of quadtree search versus walk strategies.} In this test, we compare the performance of the proposed walk strategies with that of the quadtree. We fix the Lipschitz constants of the dynamics to $L_x=L_t=2\pi$, the Courant number to $\alpha=5$, and we consider finer and finer triangular meshes with $10^5\lesssim N\lesssim 7\cdot 10^6$. The results are reported in Fig. \ref{fig:comparison}: in the left plot we show the CPU times divided by $N$, while in the right one we show the improvement factor with respect to the quadtree.

As $N$ increases, all the proposed walk strategies exhibit, as expected, an $\mathcal O(1)$ complexity, whereas CPU times for the quadtree still grow due to its $\mathcal O(\log N)$ complexity. Nevertheless, it performs better than the walking strategy $(a)$, due to the relatively large Courant number (note that the ``improvement'' factor is less than $1$). On the other hand, the improvement for the walk strategies $(b)$ and $(c)$ is apparent, with a factor between $2$ and $5$ in this setting.

\begin{figure}
	\centering
	\includegraphics[width=0.45\textwidth]{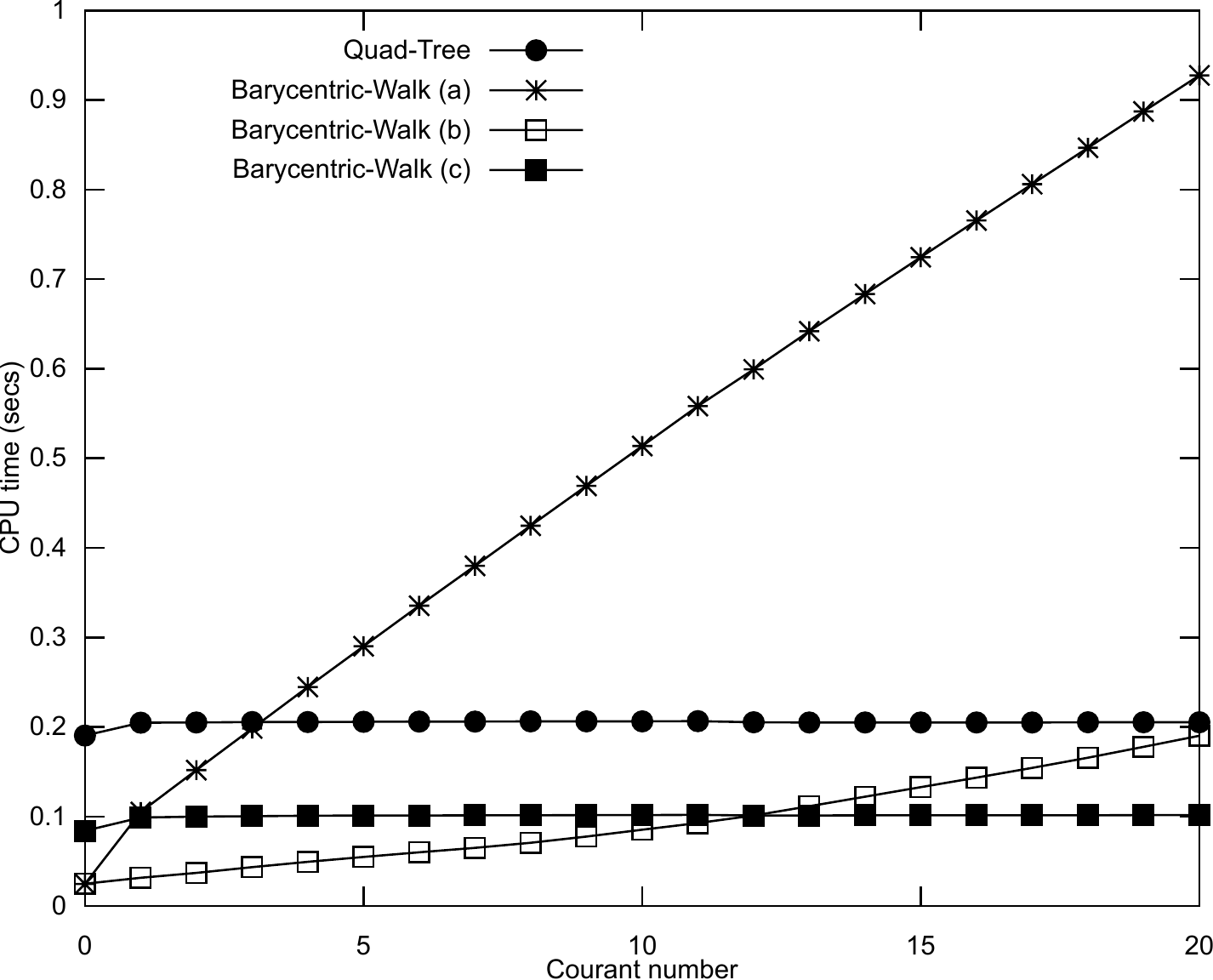} \quad
	\includegraphics[width=0.45\textwidth]{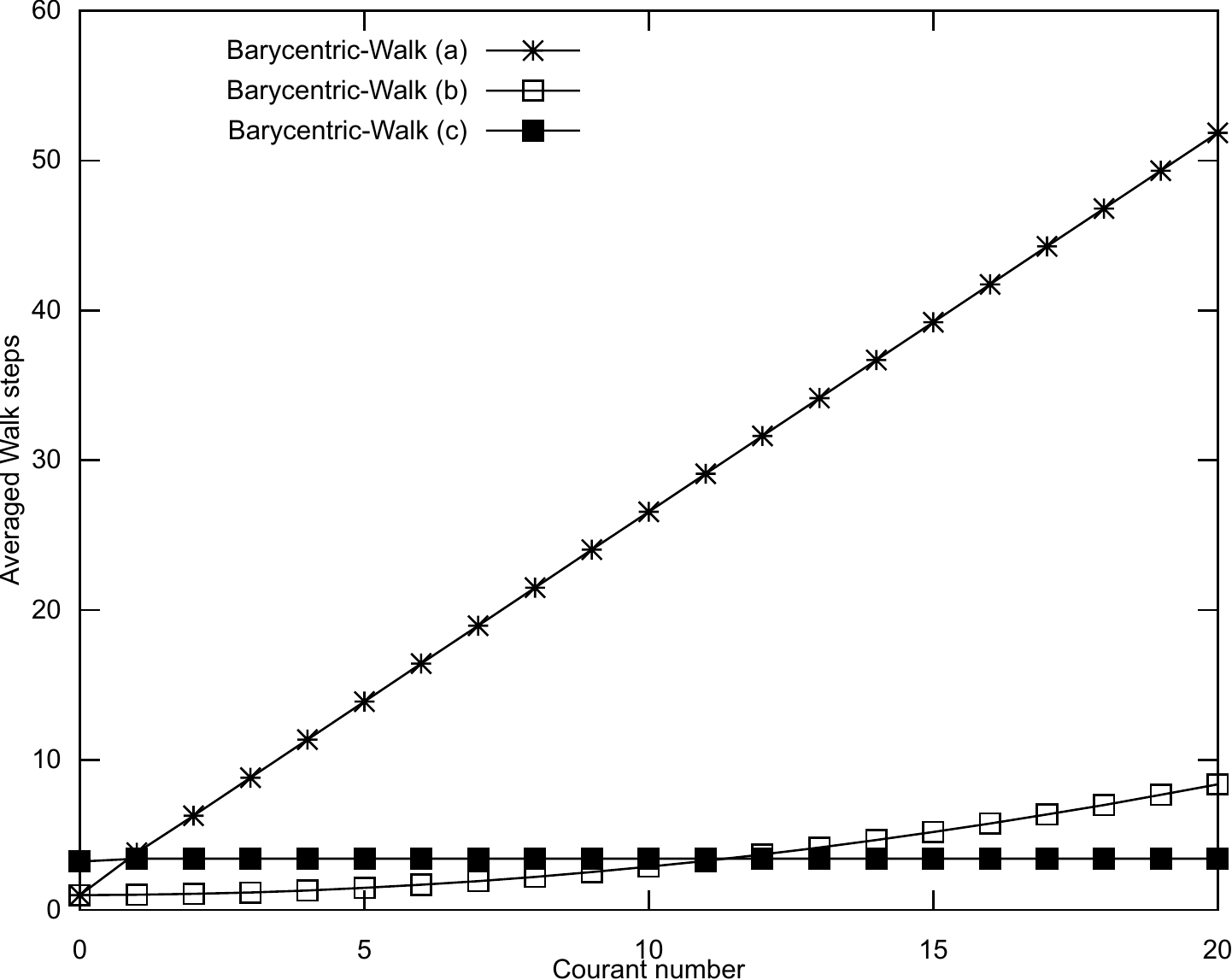}
	\includegraphics[width=0.45\textwidth]{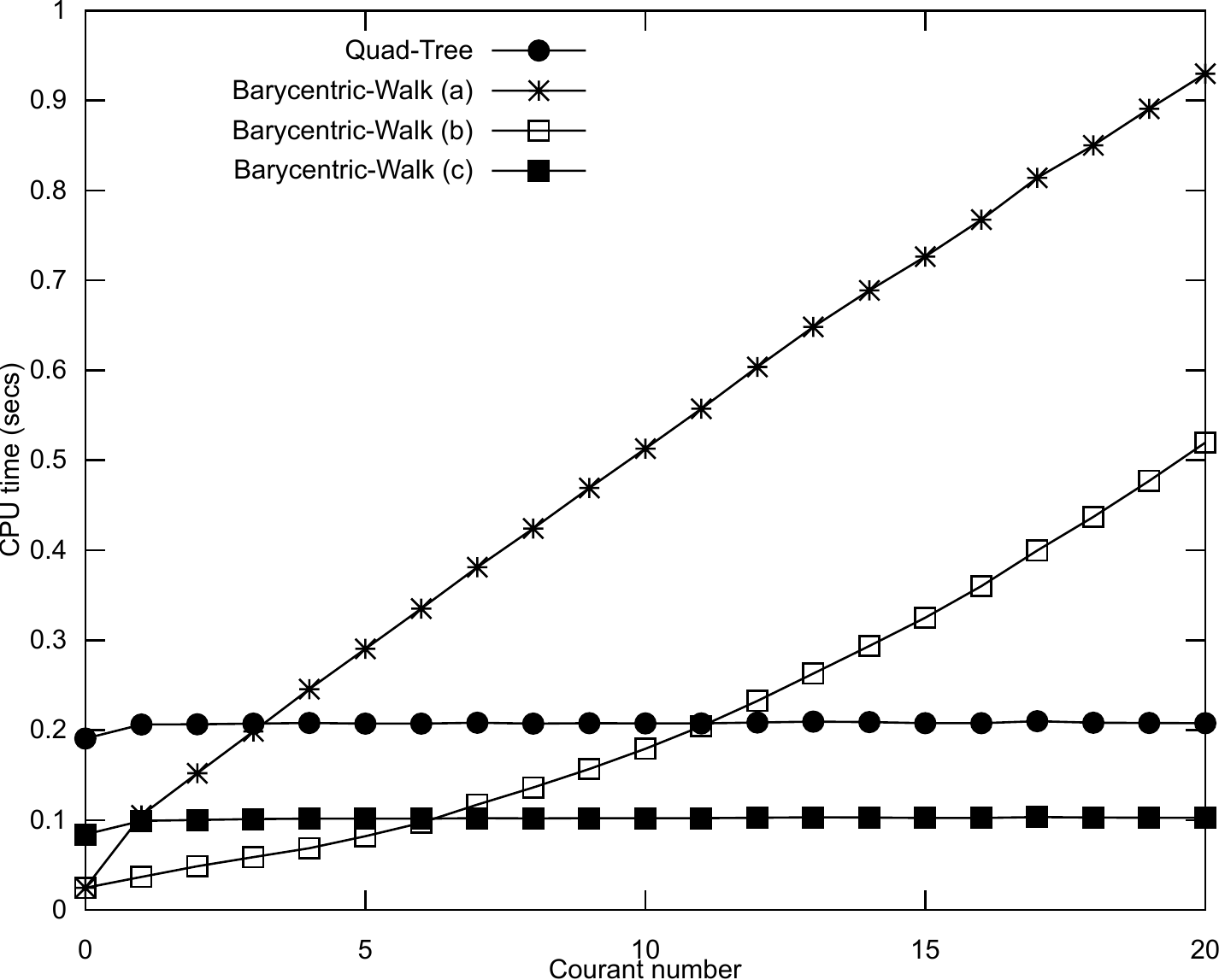} \quad
	\includegraphics[width=0.45\textwidth]{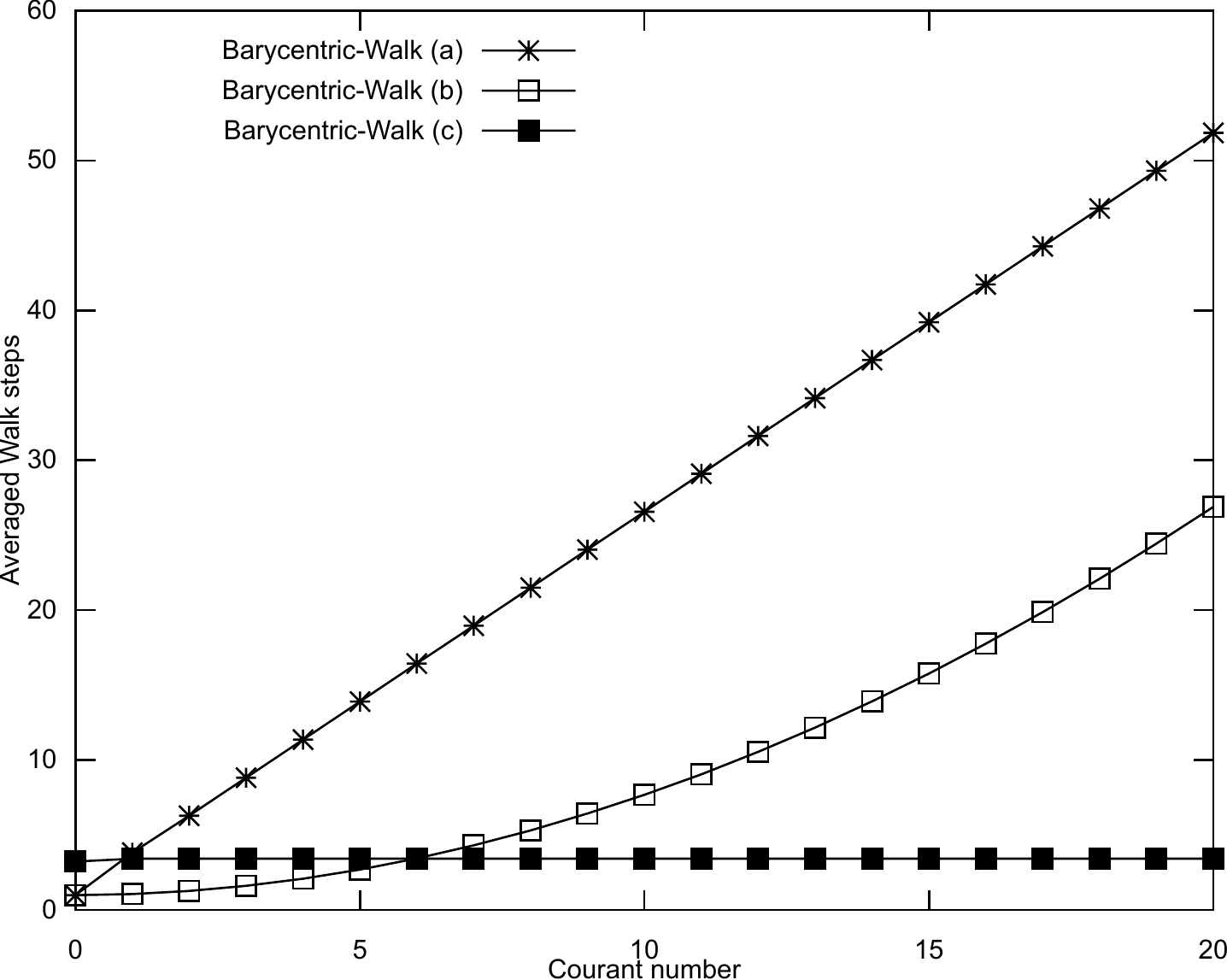}
	\caption{Search times versus Courant number (left) and averaged barycentric walk steps versus Courant number (right), for $N=1.5\cdot 10^6$, $L_x=2\pi$, and with $L_t=2\pi$ (upper plots), $L_t=8\pi$ (lower plots).}\label{fig:courant-lx1-lt4}.
	\end{figure}

\begin{figure}
	\centering
	\includegraphics[width=0.45\textwidth]{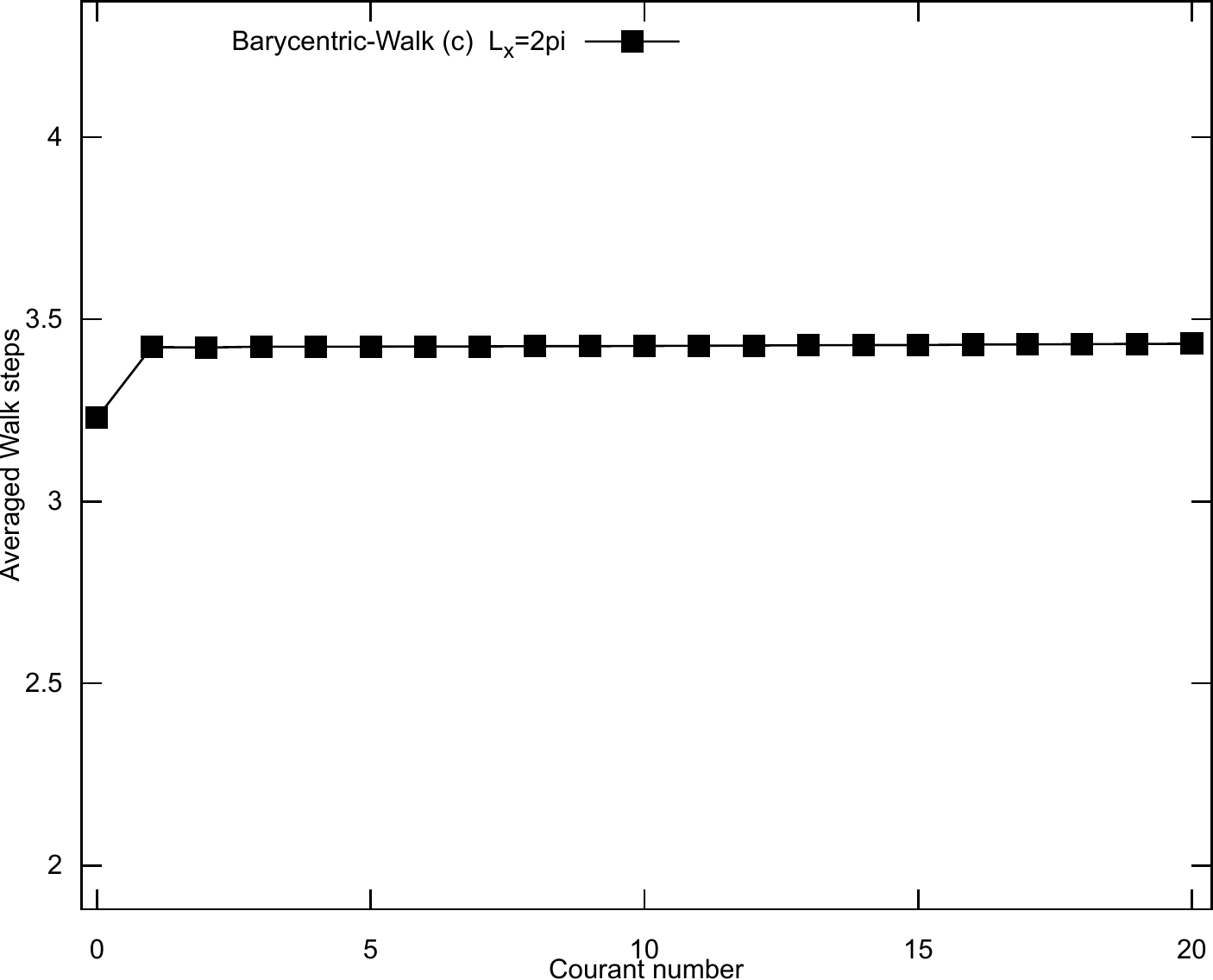} \quad
	\includegraphics[width=0.45\textwidth]{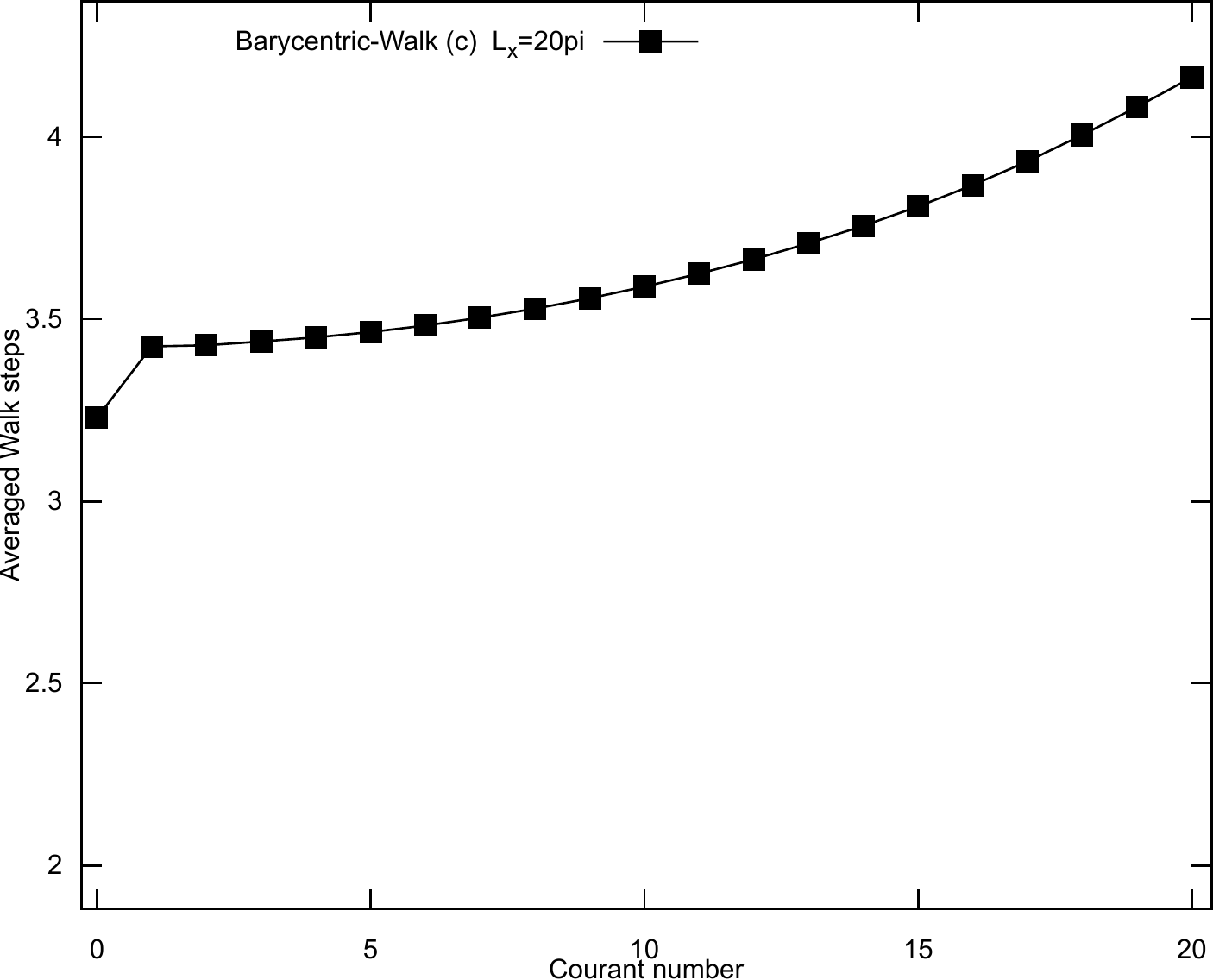}
	\caption{Averaged barycentric walk steps versus Courant number for strategy $(c)$, with $N=1.5\cdot 10^6$ and $L_x=2\pi$ (left), $L_x=20\pi$ (right).}\label{fig:courant-BWc-lx1}
	\end{figure}

\paragraph{Dependence on the Courant number.} In this test, we analyze the walk strategies in terms of the Courant number and the Lipschitz constants of the advecting dynamics. To this end, we consider a triangular mesh with a fixed number of nodes $N=1.5\cdot 10^5$ (corresponding to $A_{\max}=5\cdot 10^{-5}$ and $\Delta x=10^{-2}$), and we choose a variable Courant number in the range $0\le\alpha\le 20$. 
In this setting, the complexity for the three walk strategies can be rewritten as:
\begin{equation}\label{eq:courant-complexity}
(a)\,\,\mathcal O(C_1^B+C_2^B\alpha)\,,\qquad
(b)\,\,\mathcal O(C_1^B+C_2^BL_t\alpha^2\Delta x)\,,\qquad
(c)\,\,\mathcal O(C_1^B+C_2^B+C_2^BL_x\alpha\Delta x)\,.
\end{equation}
In the upper plots of Fig. \ref{fig:courant-lx1-lt4}, we report the results for the case $L_x=L_t=2\pi$, while in the lower plot we report the same data for $L_x=2\pi$ and $L_t=8\pi$. As expected, complexity of the quadtree search does not depend at all on the Courant number. Moreover, for $\alpha\to 0$, we recover the constant terms for the three walk strategies. In particular (see \eqref{eq:courant-complexity}), strategies $(a)$ and $(b)$ have the same constant $C_1^B$, while strategy $(c)$ shows an additional cost due to the constant $C_2^B$. 
On the other hand, as $\alpha$ increases, we clearly observe the linear behavior for strategy $(a)$, which eventually performs worse than the quadtree search. Strategies $(b)$ and $(c)$ are the most effective, due to the terms $\Delta x$ in \eqref{eq:courant-complexity}. Moreover, we recognize a quadratic behavior in $\alpha$ for strategy $(b)$, with a loss of performance at the increase of $L_t$, while the linear behavior for strategy $(c)$, in the chosen range for $\alpha$, is somewhat hidden by both its slope $\Delta x$ and the constant term. This confirms the uniform bound in \eqref{eq:uniform-bound}, since  the crossing condition $L_x\Delta t\ge 1$ for the characteristics reads, in the present case, as $2\pi\alpha\Delta x\ge 1$, namely $\alpha\gtrsim 16$.
The effect of the Lipschitz constant $L_x$ on the number of steps for the strategy $(c)$ is analyzed more in detail in Fig. \ref{fig:courant-BWc-lx1}. 
In the left plot, we use $L_x=2\pi$, and obtain an averaged number of walk steps of about $3.5$, which is in agreement with \eqref{eq:uniform-bound} in view of Remark \ref{rem:CPU_bw}. The increase in the number of steps becomes apparent when choosing a larger Lipschitz constant $L_x=20\pi$, which forces the crossing of characteristics around $\alpha\gtrsim 1.6$, so that the uniform bound on the walking steps fails, as shown in the right plot of Fig. \ref{fig:courant-BWc-lx1}. We point out, however, that this makes the scheme work in unstable (and unphysical) conditions.


In conclusion, the numerical experiments confirm that we can choose the most efficient walk strategy according to the shape of the advecting dynamics. In particular, we can conceive a hybrid (and even smarter) point location algorithm, which selects the appropriate walk strategy taking into account local information provided by the dynamics at a specific point in the domain. This direction of research is still under investigation, and we plan to address it in a forthcoming work.

%
%

\begin{figure}
\centering
\includegraphics[height=6cm]{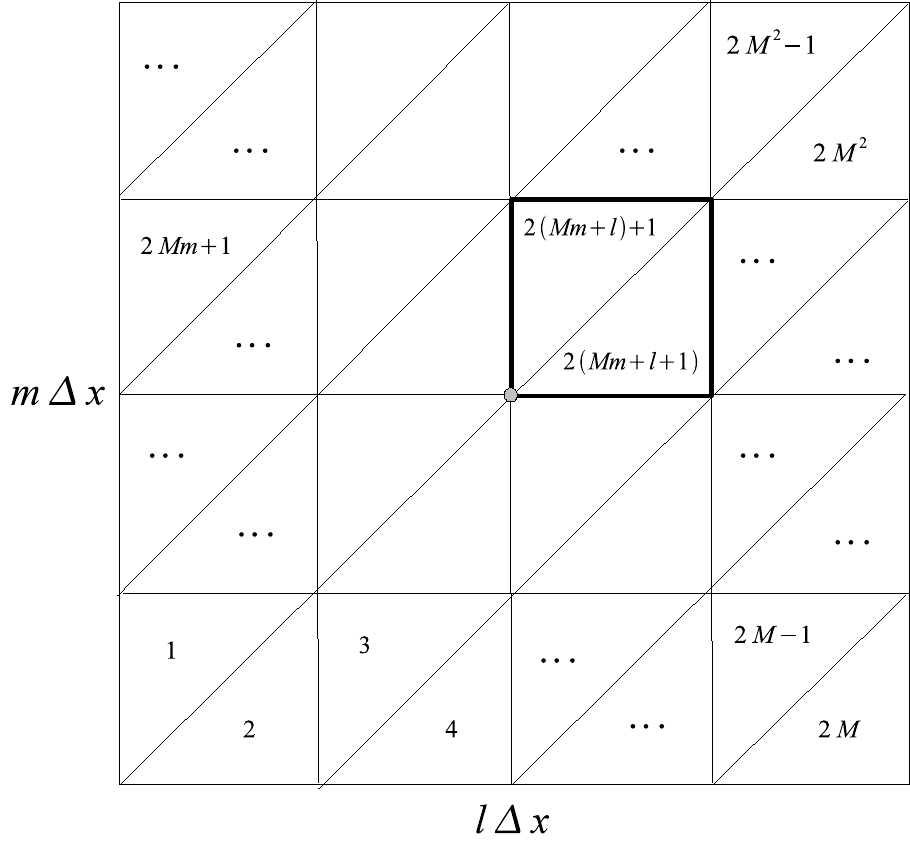}
\caption{Courant triangular grid for structured/unstructured search comparison.}\label{fig:structured-mesh}
\end{figure}

\begin{figure}
	\centering
	\includegraphics[width=0.47\textwidth]{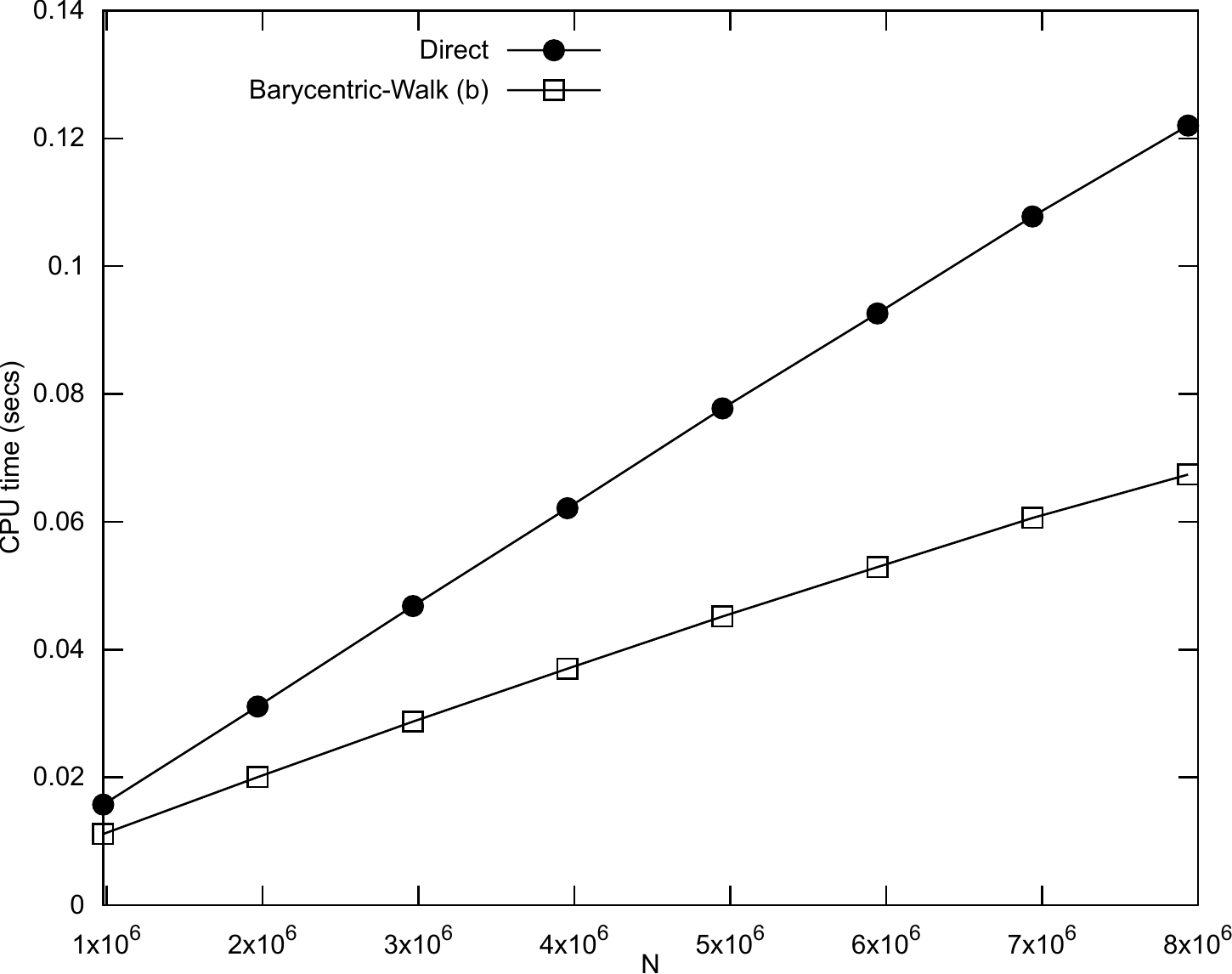} \quad
	\includegraphics[width=0.47\textwidth]{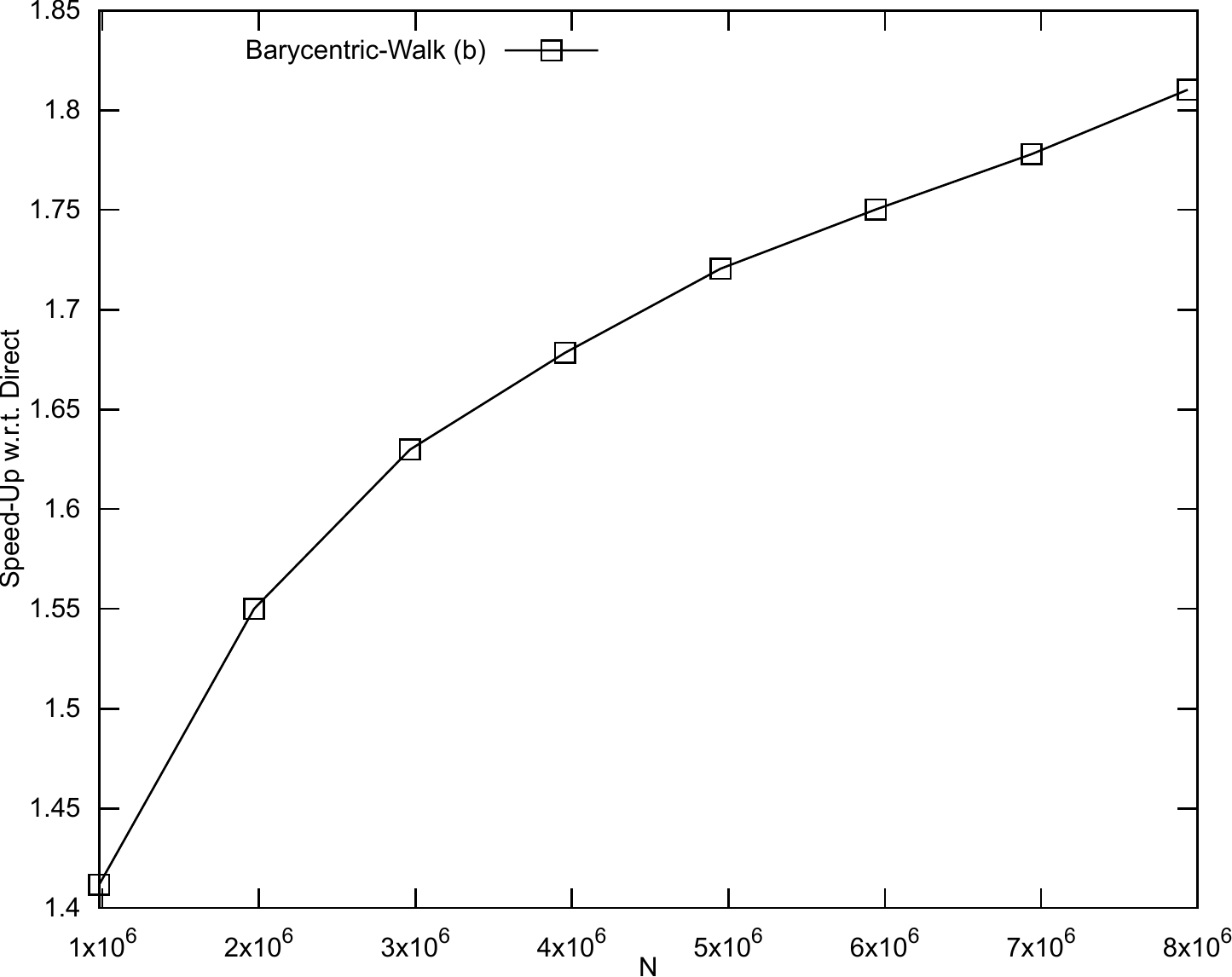}
	\caption{Search time (left) and improvement factor (right) of the barycentric walk (b) versus the direct search on a structured grid, for $10^6\lesssim N\lesssim 8\cdot 10^6$.}\label{fig:structured}
	\end{figure}

\paragraph{Comparision with direct location on a structured grid.} In this test, we compare the barycentric walk search with a direct search on a mesh which is still triangular but structured. More precisely, we consider in $[-1/2,1/2]^2$ the Courant triangulation shown in Fig. \ref{fig:structured-mesh}, with a uniform number of nodes in each dimension and a natural labelling of the corresponding triangles. 
In this setting, a given target point $X^\Delta(x_i,t_{n+1};t_n)=(\xi,\eta)$ can be directly located, with constant complexity, by the couple 
$$
l = \left\lfloor\frac{\xi+1/2}{\Delta x}\right\rfloor,\quad m = \left\lfloor\frac{\eta+1/2}{\Delta x}\right\rfloor
$$
$$
X^\Delta(x_i,t_{n+1};t_n) \in \begin{cases}
T_{2(Mm+l)+1} & \text{if } \xi+1/2-l\Delta x<\eta+1/2-m\Delta x \\
T_{2(Mm+l+1)} & \text{otherwise,} 
\end{cases}
$$
with $M$ denoting the number of elements for each side of the square.
We consider the same advecting dynamics of the previous tests, with $L_x=L_t=2\pi$ and Courant number $\alpha=5$, and we choose the walk strategy $(b)$, which achieves the best performance in this case (see again Fig. \ref{fig:courant-lx1-lt4}). For a fair comparison, we include in the direct search both the location of the triangle and the computation of the corresponding barycentric coordinates for the target point (the minimal requirement for any advecting scheme). The results are reported in Fig. \ref{fig:structured}. 
	
\noindent In the left plot we show the CPU times for the search on finer and finer meshes with a total number of nodes $10^6\lesssim N\lesssim 8\cdot 10^6$. 
As expected, both algorithms show a linear behavior in $N$, but, surprisingly, the walk strategy $(b)$ outperforms the direct search, with a factor between $1.4$ and $1.8$ (see the right plot in Fig. \ref{fig:structured}). We found out that the most expensive task in the direct search consists in the two {\tt floor} operations, and it is deeply related to the assembly code generated by the GCC compiler. A more careful test should be performed by running the code against different compilers and architectures, but this goes beyond the scope of the present paper. 
We can conclude anyway that, even in less favourable conditions, the proposed walk strategy, when implemented on unstructured meshes, has comparable performances with respect to the fully structured case, for interpolations of finite element type.


%

\paragraph{Comparison with Matlab {\tt pointLocation}.} In this test, we compare our barycentric walk search with the built-in {\tt pointLocation} Matlab function. It is known that the Matlab environment provides several facilities for practitioners, including toolboxes for generating unstructured meshes and for solving PDEs. Unfortunately, many Matlab functions (as the general-purpose {\tt pointLocation}) are closed-source, pre-compiled, and they also lack documentation on algorithmic details, hence there is no chance to modify them for specific tasks. Here, we provide some hints to implement our barycentric walk algorithm with few lines of code in Matlab, then we evaluate its performances versus the {\tt pointLocation} function. To this end, we employ the Matlab command {\tt mex} to build a MEX function from our C++ implementation, namely a binary file that can be called, as any Matlab built-in function, by a Matlab script. Starting from the Matlab {\tt triangulation} data structure, containing point coordinates of the nodes and vertex indices of the triangles in the mesh, we add three additional fields: a list containing the indices of the initial triangles for the barycentric walk (one index per node), a list containing the index triplets of neighboring triangles (one triplet per triangle), and a list containing the indices of the parent nodes for the walk strategy $(c)$ (one index per node). In particular, the first list can be constructed choosing a random triangle from the output of the Matlab function {\tt vertexAttachments}, the second list is simply the output of the Matlab function {\tt neighbors}, while the third list can be obtained, starting from a root node, using {\tt vertexAttachments} to find recursively the first neighbors of nodes already inserted in the list (some care must be taken to avoid duplicates). Then, we design our MEX function  {\tt pointLocationBW} with a syntax similar to {\tt pointLocation}:
$$
{\tt [I,B] = pointLocationBW(TBW,Q,s)}\,,
$$
where {\tt TBW} is the extended triangulation data structure, {\tt Q} the list of query points, {\tt s} the chosen walk strategy ({\tt 'a'},{\tt 'b'} or {\tt 'c'}), while {\tt I} is the output list of the triangles enclosing the query points, {\tt B} the corresponding list of barycentric coordinates. 
This function is still in beta version, but available for the interested readers on reasonable request.
     
\begin{figure}
	\centering
	\includegraphics[width=0.47\textwidth]{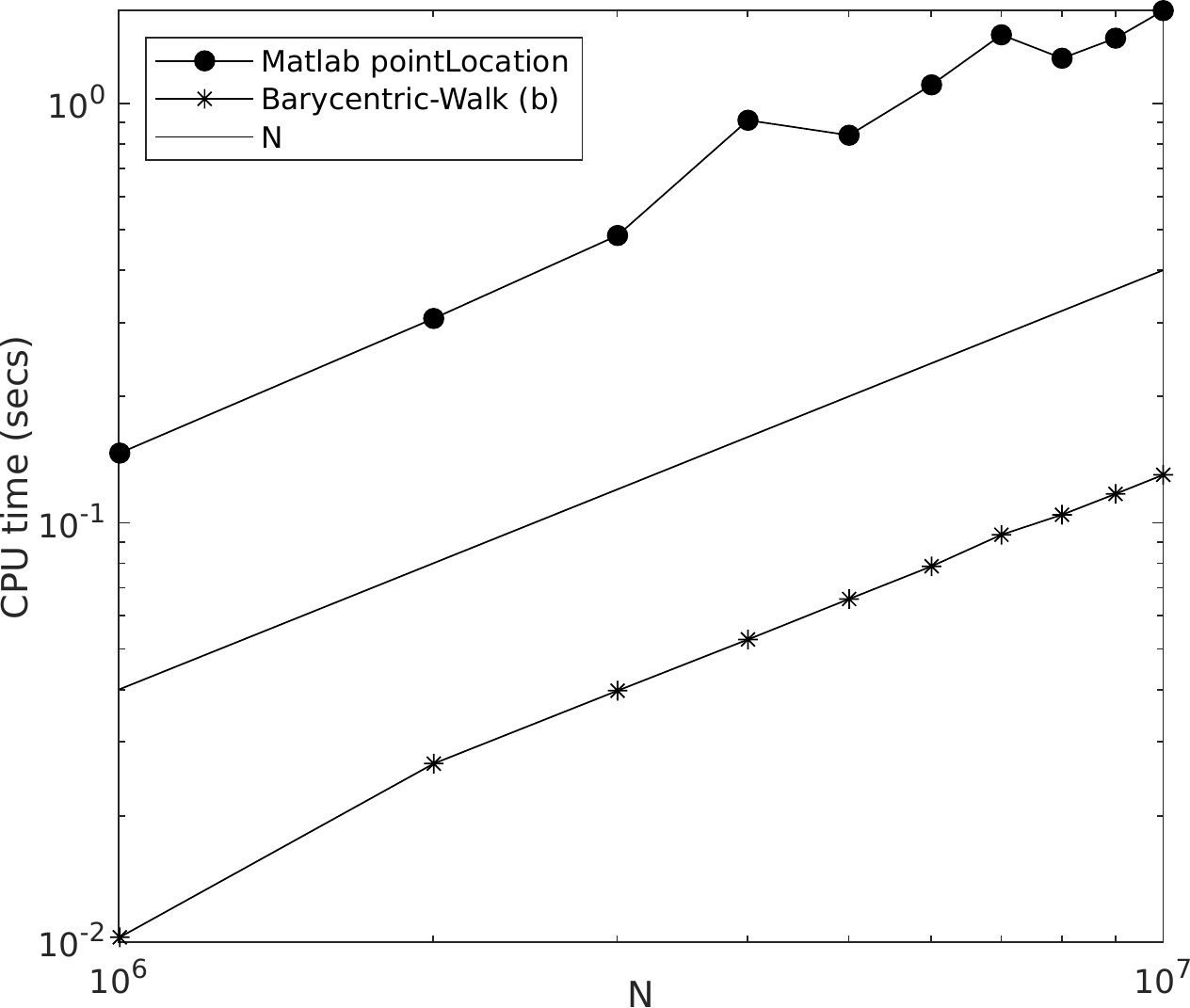} \quad
	\includegraphics[width=0.47\textwidth]{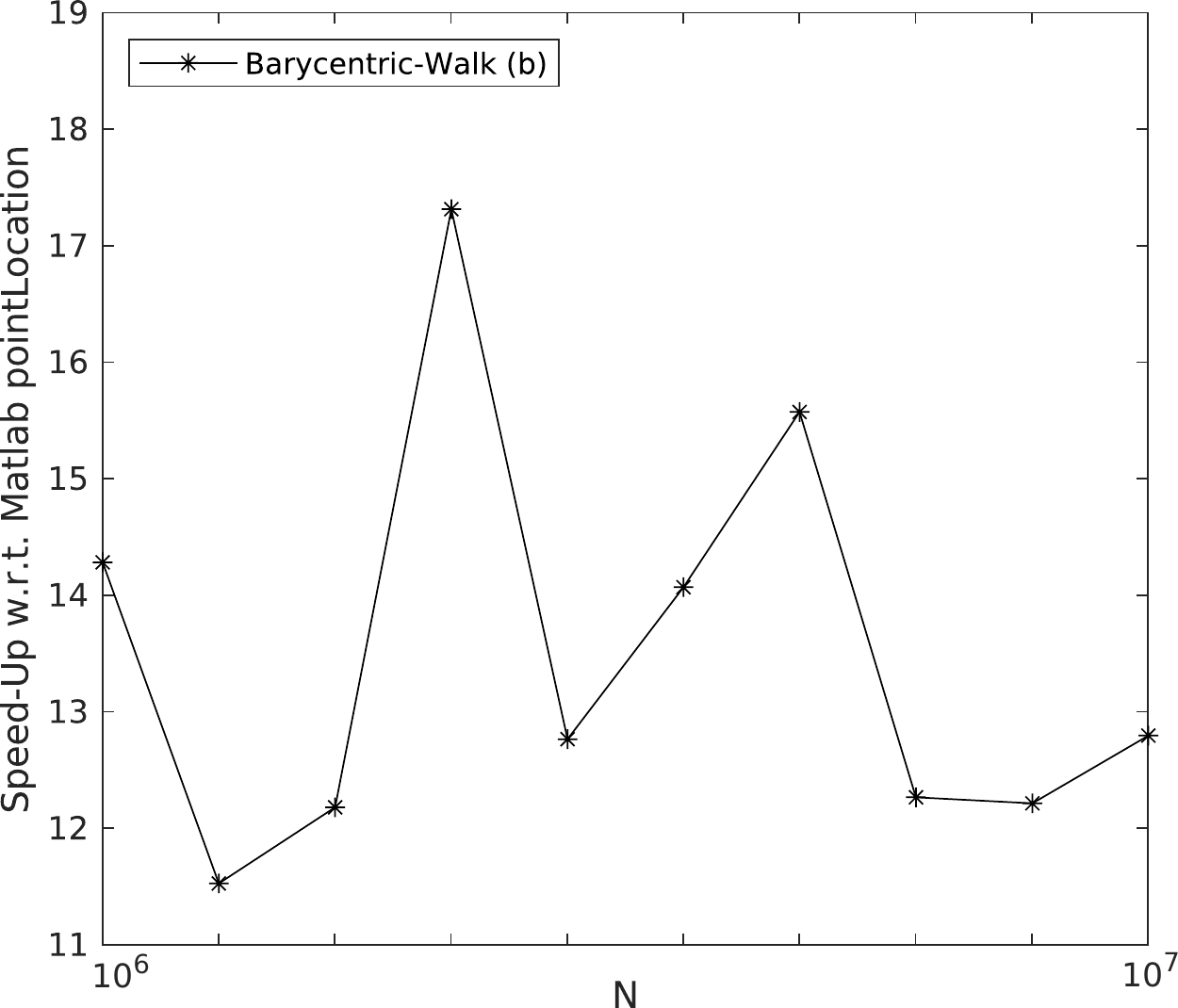}
	\caption{Search time (left) and improvement factor (right) of the barycentric walk (b) versus the Matlab {\tt pointLocation}, for $10^6\lesssim N\lesssim 10^7$.}\label{fig:matlab}
	\end{figure}

Now, we set the same parameters of the previous test, namely $L_x=L_t=2\pi$, $\alpha=5$,
then we choose the walk strategy $(b)$ 
and run the code on Matlab version {\tt R2021a}.  
The results are reported in Fig. \ref{fig:matlab}. 

\noindent In the left plot we show the CPU times for the search on finer and finer meshes with a total number of nodes $10^6\lesssim N\lesssim 10^7$. We observe that also the Matlab {\tt pointLocation} seems to have a linear complexity in $N$, and this would suggest that its black-box algorithm might not be based on a quadtree structure. The function {\tt pointLocationBW} improves the {\tt pointLocation} CPU times by a factor ranging from $11$ to $18$, as shown in the right plot.

\begin{figure}
	\centering
	\includegraphics[width=0.47\textwidth]{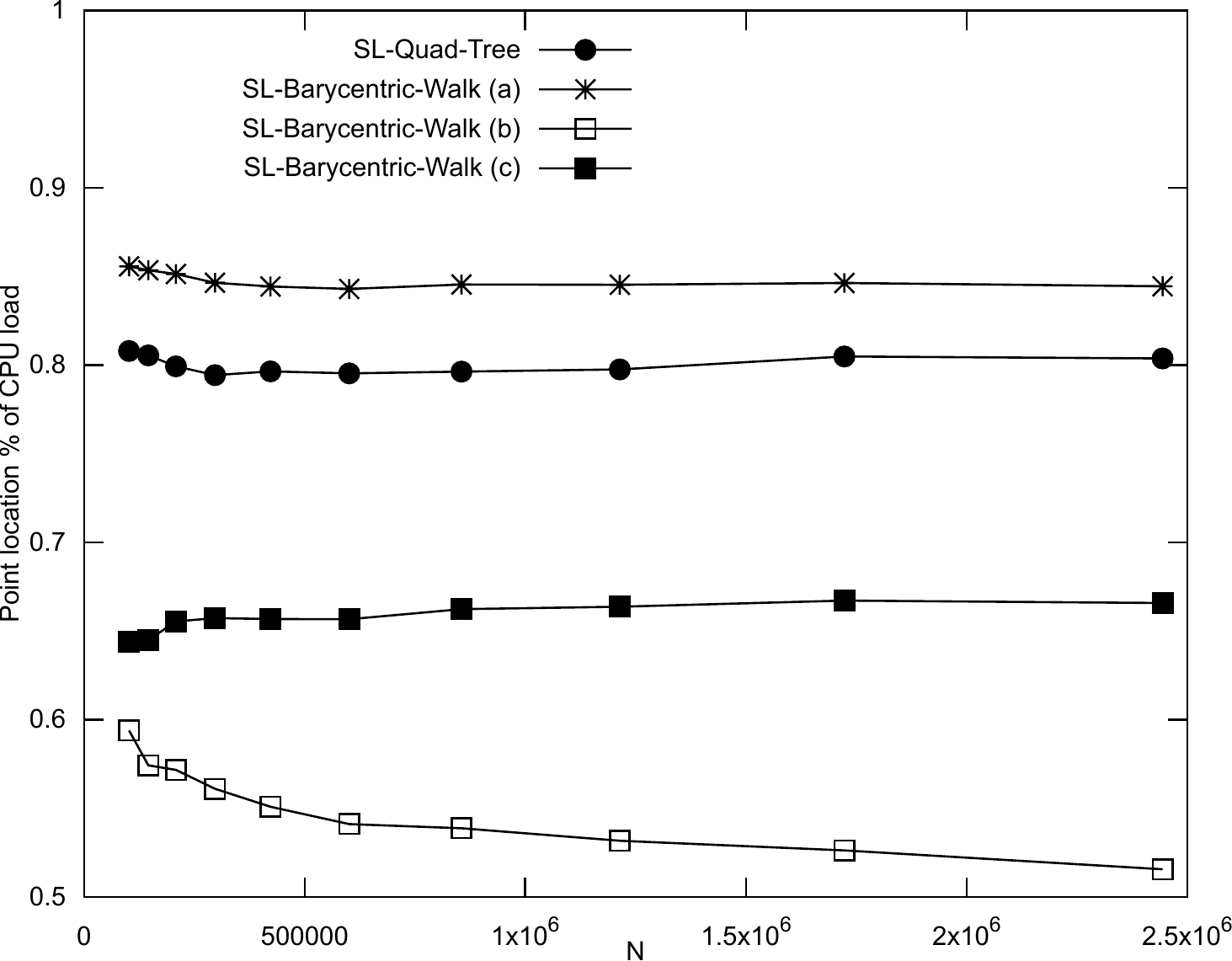}
	\caption{Percentage of CPU load due to point location, for different algorithms in a SL scheme, versus number of grid nodes $10^5\lesssim N\lesssim 2.5\cdot 10^6$.}\label{fig:SL-times}
	\end{figure}

\paragraph{Point location CPU load in the SL scheme.}
In this last test, we combine the point location provided by the quadtree search and the proposed walk strategies, with the SL scheme ($\mathbb P_1$ interpolation + Euler tracking of characteristics) illustrated in Algorithm \ref{ALG2}. The aim is to measure, for the different algorithms, which percentage of the total computational load is due to the point location.  
To this end, we consider the same parameters for the advecting dynamics of the previous tests ($L_x=L_t=2\pi$, $\alpha=5$), and we choose a Gaussian-like initial datum $u_0$. The results for different meshes of size $10^5\lesssim N\lesssim 2.5\cdot 10^6$ are reported in Fig. \ref{fig:SL-times}.  
We remark that here the total CPU time for each run includes,  on the whole mesh and for all the time steps, the computation of the query points $\mathcal Q^n$, the point location and the interpolation of the solution. Then we show the ratio between the point location time and the total CPU time. We observe that the quadtree point location achieves about $80\%$ of CPU load, against the $85\%$ of the walk strategy $(a)$. On the other hand, we get about $65\%$ for strategy $(c)$, while for strategy $(b)$ the percentage drops between $50\%$ and $60\%$ (in particular it decreases as $N$ increases before saturating, whereas the CPU load is actually constant for the other algorithms). This is not surprising, since strategy $(b)$ is designed to take advantage from small variations of the dynamics with respect to time. Note that in this test we have a moderate value for $L_t$, while the Courant number and the final time are kept fixed, so that the number of time steps increases with $N$. This implies that most characteristics eventually  fall in the same triangle for more and more time steps.

\section{Conclusions}\label{sec:conclu}

In this paper, we have analyzed in detail the complexity issues related to characteristics location in SL-type schemes on 2-D unstructured triangular space grids. We have proposed two new and clever choices of the initial element for the barycentric walk point location, which may bring this algorithm to a higher degree of efficiency with respect to the recipes typically used so far, in particular when the advection term has a slow variation with respect to the space and/or time variable.

Although the analysis has been carried out in a specific setting (triangular grids, two-dimensional problems), it is not difficult to extend the technique to more general situations, in particular to Voronoi meshes, as well as to three-dimensional problems. The choice of a walk algorithm different from the barycentric walk is also possible (see \cite{DPT01}), especially to treat the case of less regular space grids than the ones we have used here.

\subsection*{Acknowledgements}

This work has been partially supported by the PRIN 2017 project {\it ``Innovative Numerical Methods for Evolutionary Partial Differential Equations and Applications''}, by the INdAM--GNCS project {\it ``Approssimazione numerica di problemi di natura iperbolica ed applicazioni''} and by Roma Tre University.

\noindent We thank Dr.\,Beatrice Beco and Dr.\,Lorenzo Della Cioppa for taking part in the first steps of this work.

\end{document}